\documentclass[AMA,STIX1COL,doublespace]{WileyNJD-v2}
\articletype{Research Article}

\usepackage{amsmath,amsfonts,amsthm}
\usepackage{graphicx}
\usepackage{subcaption,fancybox}
\usepackage{bm}
\usepackage{verbatim}
\usepackage{cancel}
\usepackage{color,colordvi}
\usepackage{tikz}
\usetikzlibrary{plotmarks}
\newcommand\marksymbol[2]{\tikz[#2,scale=1.2]\pgfuseplotmark{#1};}

\newcommand{\suku}[1]{{#1}} 
\newcommand{\proofcor}[1]{{#1}} 

\newcommand{\fref}[1]{Fig.~\ref{#1}}
\newcommand{\tref}[1]{Table~\ref{#1}}
\newcommand{\sref}[1]{Section~\ref{#1}}
\newcommand{\vm}[1]{\bm{#1}}
\newcommand{\vx}{\vm{x}}
\newcommand{\face}{\texttt{f}}

\renewcommand{\Re}{{\mathbb{R}}}
\newcommand{\xivec}{{ \vm{\Xi} }}
\newcommand{\etavec}{{ \vm{H} }}
\newcommand{\zetavec}{{ \vm{Z} }}

\graphicspath{ {./figs/} }

\begin{document}
	
\title{Virtual elements on agglomerated finite elements to increase 
       the critical time step in elastodynamic simulations}
	
\author[1]{N. Sukumar$\mbox{}^{*,}$}
\author[2]{Michael R. Tupek}
	
\authormark{N. SUKUMAR \sc{and} MICHAEL R. TUPEK}
	
\address[1]{\orgdiv{Department of Civil and Environmental Engineering},
	\orgname{University of California},
	\orgaddress{Davis, CA 95616, USA}}
\address[2]{\orgname{Sandia National Laboratories},
	\orgaddress{Albuquerque, NM 87185}, \country{USA}}
\corres{$\mbox{}^*$N. Sukumar, Department of Civil and 
	Environmental Engineering,
	University of California, One Shields Avenue, Davis, 
	CA 95616, USA \\
	\email{nsukumar@ucdavis.edu}}

\abstract{
In this paper, we use the first-order virtual element method (VEM) 
to investigate the effect of shape quality of polyhedra in the estimation of the critical time step for explicit three-dimensional elastodynamic finite element (FE) simulations. Low-quality finite elements are common when meshing realistic complex components, and while tetrahedral meshing technology is generally robust, meshing algorithms cannot guarantee high-quality meshes for arbitrary geometries or for non-water-tight computer-aided design models. For reliable simulations on such meshes, we consider FE meshes with tetrahedral and prismatic elements that have badly-shaped elements---tetrahedra with dihedral angles close to $0^\circ$ and $180^\circ$, and slender prisms with triangular faces that have short edges---and agglomerate such `bad' elements with neighboring elements to form a larger polyhedral virtual element.  On each element,  
the element-eigenvalue inequality is used to estimate the critical time step.  For a suite of illustrative finite element meshes with $\epsilon$ being a mesh-coordinate parameter that leads to poor mesh quality, we show that adopting VEM on the agglomerated polyhedra yield critical time steps that are insensitive as $\epsilon \rightarrow 0$.
\suku{
The significant reduction in solution time on meshes with agglomerated virtual elements vis-{\`a}-vis tetrahedral meshes  is demonstrated 
through explicit dynamics simulations on a tapered beam.}
}

\maketitle

\section{Introduction}\label{sec:intro}
Over the past fifteen years there has been a sustained effort in the development of new finite element formulations on meshes consisting of arbitrary polytopal (polygonal and polyhedral) elements.\cite{Hormann:2017:GBC}
The virtual element \suku{method} (VEM)~\cite{Beirao:2013:BPV} is a recent contribution in higher order stable Galerkin discretizations on polytopal meshes to solve boundary-value problems. It has its roots in mimetic finite-difference schemes that are based on an algebraic construction for the system (stiffness, mass) matrices.\cite{Beirao:2014:MFD} The VEM endows a variational framework for the first-order mimetic scheme, and in doing so, it provides a route to many new formulations.
In the VEM, the basis functions are defined as the solution of a local
elliptic partial differential equation, and are never explicitly
computed in the implementation of the method. 
For this reason, they are referred to as \textit{virtual}, and the finite
element space of the VEM as the \emph{virtual element space}.
As in the finite element method (FEM), the local discretization space
$V_k(E)$ ($k$ is the order of the element $E$) from each element is glued
together to form the conforming global space $V^h$. However, distinct
from the standard FEM, the trial and test functions in $V_k(E)$
consist of polynomials of order less than or equal to $k$ and of 
additional nonpolynomial functions.
Since the virtual basis functions are unknown in $V_k(E)$, the VEM uses their elliptic polynomial projections to build the bilinear form (stiffness
and mass matrices) and the continuous linear functional (forcing terms) of the variational formulation.  For problems such as the Poisson equation or
linear elastostatics, such projections are computable from the degrees of freedom without introducing any further approximation error and are used 
to decompose the bilinear form on each element into two parts: the consistent term that approximates the stiffness matrix on a given polynomial space and the correction term that ensures coercivity (stability with the proper scaling). As in the FEM, element-level assembly
procedures are used to obtain the global system matrices. For a simplex in $\Re^d$, the first-order~($k = 1$) VEM is identical to linear FEM.  The VEM can be viewed as a stabilized hourglass 
control finite element method~\cite{Flanagan:1981:AUS} on convex and 
nonconvex polytopes.\cite{Cangiani:2015:HSV}

In comparison to the standard finite element method, 
the simplicity and flexibility of the VEM has attracted attention from both applied mathematicians and engineers who have furthered the theory and applications of the method. The
ability to devise accurate and stable numerical simulations
on polygonal and polyhedral meshes with arbitrarily-shaped convex 
and nonconvex elements has led to many new applications of the method.
Initially, the emphasis of the virtual element method was on low- and high-order 
formulations for two- and three-dimensional scalar elliptic problems.\cite{Ahmad:2013:EPV,Beirao:2014:HGV,Beirao:2017:HOV,Dassi:2018:EHO}
Over the past five years, greater emphasis has been placed on vectorial problems such as fluid flow, and linear and nonlinear deformation of solids. Among these, 
we mention
applications of the VEM in 
linear elastostatics,~\cite{Beirao:2013:VEL,Gain:2014:VEM,Artioli:2017a:AO2,Mengolini:2019:EPV} linear  elastodynamics,\cite{Park:2019:ONM,Park:2020:NRE,Antonietti:2021:AVE} and finite strain elastodynamics\cite{Cihan:2020:VEF}
that are pertinent to this contribution. 
Park et al~\cite{Park:2019:ONM,Park:2020:NRE} were the first to
show that the first-order 
VEM with explicit
time integration is accurate and converges optimally for linear elastodynamic simulations 
over convex and nonconvex polygonal and polyhedral meshes.
Antonietti et al.\cite{Antonietti:2021:AVE} proved stability and convergence of
the semi-discrete arbitrary-order virtual element
method for two-dimensional linear elastodynamic problems, and Cihan et al.\cite{Cihan:2020:VEF} extended the
low-order virtual element method for two- and three-dimensional 
finite strain elastodynamic simulations.

Low-quality finite elements are common when meshing realistic complex 
components, and while tetrahedral meshing technology is generally robust, meshing algorithms cannot guarantee high-quality meshes for arbitrary geometries or for non-water-tight computer-aided design models.\cite{Koester:2019:ADS}
In explicit linear elastodynamic computations,
the presence of tetrahedral slivers (dihedral
angles close to $0^\circ$ and $180^\circ$) in a finite element mesh
dramatically decreases the critical time step, making simulation progress practically impossible. 
In this paper, we use the first-order VEM to
 investigate the effect of shape quality of polyhedra (formed via
agglomeration of poor-quality elements in tetrahedral-dominant finite element meshes) 
in the estimation of the critical time step for explicit 
dynamics simulations. Agglomeration of tetrahedral elements to form
polyhedra have been introduced for 
discontinuous Galerkin simulations,~\cite{Bassi:2012:OTF,Cangiani:2014:HPD}
and more recently in polyhedral finite element computations over tetrahedral
meshes.\cite{Bishop:2020:PFE}
The structure of the remainder of this paper
follows. The elastodynamic
eigenproblem (strong and weak forms) is presented in~\sref{sec:eigproblem}. In Section~\ref{sec:math}, we present some of the essentials on polynomial basis sets and projection operators. Then, in Section~\ref{sec:vem}, the
main concepts in the VEM---definition of the virtual element space and degrees of freedom---are presented. The discrete generalized eigenproblem is also stated. In~\sref{sec:implementation}, the numerical implementation of three-dimensional VEM for linear solid mechanics problems is presented, including expressions for the element stiffness and element mass matrices, and details on the numerical integration using the homogeneous numerical integration scheme.  It is known that triangles with an internal angle close to $180^\circ$ 
(flat triangles) adversely affect interpolation accuracy, whereas an internal
 angle close to $0^\circ$ (and the others not large) worsen stiffness-matrix conditioning.\cite{Shewchuk:WGL:2002}  The dihedral angle in a tetrahedron dictates element quality from
the viewpoints of interpolation and stiffness-matrix conditioning.\cite{Klingner:2008:ITM} 
Establishing such shape quality measures for polygons and polyhedra is challenging.\cite{Gillette:2017:SQG,Attene:2021:BPQ}
In the numerical results presented in~\sref{sec:results}, we consider polyhedral elements that are formed by agglomerating  bad-quality tetrahedral elements (small edges, small faces, dihedral angles close to $0^\circ$ and $180^\circ$) with neighboring tetrahedra of good quality. We refer to such an element as a \emph{polyhedral virtual element}.
The elastodynamic eigenproblem on each element is solved, and the element eigenvalue inequality~\cite{Fried:1972:BEE,Lin:EET:1989,Lin:BEF:1991} is used to obtain a lower bound estimate of the critical time step.  Comparisons of the critical time step on a suite of bad-quality tetrahedral finite elements are made with a polyhedral virtual
element to assess the performance of the VEM. 
\suku{Finally, we present 
explicit dynamics simulations on a three-dimensional 
tapered beam to demonstrate that the
increase in critical time step using
polyhedral virtual elements leads to significantly faster computer simulations when compared to 
tetrahedral finite elements.}
We close with a summary of our main findings and conclusions in Section~\ref{sec:conclusions}. 

\section{Strong and Weak Forms of the Three-Dimensional Linear 
         Elastodynamic Eigenproblem}\label{sec:eigproblem}
The motion of waves in a linearly elastic solid is governed by the equations of elastodynamics.  The governing field equations in a domain $\Omega \subset \Re^3$ are:

\begin{subequations}\label{eq:2d-governing}
\begin{align}
  \nabla \cdot \vm{\sigma} & = \rho \ddot{\vm{u}}, \ \ \vm{\sigma} =
  \vm{\sigma}^T
    & \text{(balance of linear and angular momentum)} , \\
  \vm{\sigma} & = \mathbb{C} : \vm{\varepsilon} 
    & \text{(linear elastic constitutive law)} , \\
  \vm{\varepsilon} & = \frac{1}{2} \left( \nabla \vm{u} 
      + ( \nabla \vm{u} )^T \right)
    & \text{(small-strain kinematics)} ,
\end{align}
\end{subequations}
where $\vm{\sigma} := \vm{\sigma} (\vx, t)$ is the Cauchy stress,
$\rho := \rho (\vx)$ is the density, $\vm{u} := \vm{u} (\vx, t)$ is the displacement field,
$\vm{\varepsilon} := \vm{\varepsilon} (\vx, t)$ is the strain field, and $\mathbb{C}$ is the linear elastic material moduli tensor. Assuming
time-harmonic fields with angular frequency $\omega$, the field quantities can be decomposed to yield
\begin{equation}\label{eq:time-harmonic-fields}
  \vm{u} (\vx, t) = \vm{u} (\vx) \exp (-i \omega t) , \quad 
  \vm{\sigma} (\vx, t) = \vm{\sigma} (\vx) \exp (-i \omega t) , \quad 
  \vm{\varepsilon} (\vx, t) = \vm{\varepsilon} (\vx) \exp (-i \omega t) .
\end{equation}
Combining~\eqref{eq:2d-governing} and~\eqref{eq:time-harmonic-fields}, we
obtain the elastodynamic eigenproblem for free vibrations:
\begin{subequations}\label{eq:strong}
\suku{
\begin{align} 
  \nabla \cdot \vm{\sigma}(\vx) &= -\lambda \rho(\vx) \vm{u}(\vx) \ \ \textrm{in } \Omega, \\
\intertext{where $\lambda = \omega^2$ is the eigenvalue, and homogeneous 
           displacement boundary conditions are imposed on the boundary:}
  \vm{u} &= \vm{0} \ \ \textrm{on } \partial \Omega .
\end{align}
}
\end{subequations}

The weak form of the eigenproblem in~\eqref{eq:strong} is: 
find the displacement field (trial \suku{eigenfunction})
 \suku{$\vm{u} \in V = [H_0^1(\Omega)]^3$} and eigenvalue $\lambda = \omega^2 \in
\Re_+$, such that
\begin{subequations}\label{eq:weak}
\begin{align}\label{eq:weak-a}
  a(\vm{u}, \vm{v}) &= \lambda \,
    b(\vm{u}, \vm{v}) \quad
    \forall \vm{v} \in V,
\intertext{where $\vm{v}$ is the test displacement field,
\suku{$H_0^1(\Omega)$ is the Hilbert space that consists of
functions in $\Omega$ that have
square integrable derivatives up to order 1 and which vanish on
$\partial \Omega$}, and the
bilinear forms $a(\cdot,\cdot)$ and $b(\cdot,\cdot)$ are given by} 
  \label{eq:weak-b}
  a ( \vm{u}, \vm{v} ) := \int_\Omega \vm{\sigma} (\vm{u}) :
  \vm{\varepsilon} (\vm{v}) & \, d \vx , \quad
  b ( \vm{u}, \vm{v} ) := \int_\Omega \rho \, \vm{u} 
     \cdot \vm{v} \, d \vx .
\end{align}
\end{subequations}

\section{Mathematical Preliminaries}\label{sec:math}
\subsection{Polynomial basis}
Let $\Omega \subset \Re^d$ ($d = 2,3$) be the problem domain, and
${\cal T}^h$ a decomposition of $\Omega$ into nonoverlapping \emph{ simple} 
polygons or polyhedra. The number of nodes in ${\cal T}^h$ is $N$. 
This definition for the elements permits the inclusion of elements with
consecutive coplanar edges and faces.
We refer to $E \in {\cal T}^h$
as an \emph{element}. Let $|E|$ be the area in 2D and volume 
in 3D of $E$. The diameter of $E$
is denoted by $h_E$ and $\vx_E$ is the centroid (barycenter) of $E$.
The vertices of $E$ are denoted by $v_i$, and the coordinate of vertex 
$v_i$ by $\vx_i \equiv (x_i,y_i)$ in 2D and $\vx_i \equiv  (x_i,y_i,z_i)$ in 3D.
A polygon $E$ has $N_E$ vertices and $N_E$ edges, with the edges
    denoted 
by $e_i$ ($i=1,2,\ldots,N_E$). 
In 3D, a polyhedron $E$ has $N_E$ vertices, and 
is described by its faces (each is assumed to be a planar polygon) that are labeled as $\face_i$ $(i = 1,2,\ldots,N_E^\face)$,
where $N_E^\face$ is the number of faces of the polyhedron.  
Each
face $\face_i$ is given by a vertex connectivity in
counter-clockwise orientation such that the normal to the face points out of $E$.

Let $\mathbb{P}_k (E)$ be the function space on $E$ that consists of 
all polynomials of order less than or equal to $k$. By convention, $\mathbb{P}_{-1} = \{0\}$.  
The dimension of $\mathbb{P}_k(E)$ is denoted by $\textrm{dim}\, \mathbb{P}_k(E)$. In one dimension,  
$\textrm{dim} \, \mathbb{P}_k(E) = k+1$; in two dimension, 
$\textrm{dim} \, \mathbb{P}_k(E) = (k+1)(k+2)/2$, and in three dimension, 
$\textrm{dim} \, \mathbb{P}_k(E) = (k+1)(k+2)(k+3)/6$.
The set consisting of the scaled monomials of order less than 
or equal to $k$ on $E$ is defined as ${\cal M}_k(E)$. For example, in
two dimensions
\begin{equation*}
{\cal M}_1(E) = \left\{ 1 , \ \frac{x - x_E}{h_E} ,  \dfrac{y - y_E}{h_E} 
 \right\}, 
\end{equation*}
and in three dimension
\begin{equation*}
{\cal M}_1(E) = \left\{ 1 , \ \frac{x - x_E}{h_E} ,  \dfrac{y - y_E}{h_E} ,  \frac{z-z_E}{h_E}
                 \right\}
\end{equation*}
are scaled monomial basis sets of order 1. 
Note that ${\cal M}_1(E)$ is a basis for $\mathbb{P}_1(E)$, and all elements of ${\cal M}_1(E)$ are ${\cal O}(1)$. The scaled monomial basis set is used in the virtual element computations for scalar 
problems.  In the literature, $m_\alpha$ is used to refer to an 
element of ${\cal M}_k(E)$, and 
${\cal M}_k^* (E)$ is used to denote the basis set on $E$ that contains all monomials that are exactly of order $k$.  

Let $\vm{P}_1(E) = \mathbb{P}_1(E) \otimes \vm{I}_d$
($\vm{I}_d$ is the $d \times $d identity matrix and
$\otimes$ is the Kronecker product) be the
\suku{vector polynomial basis defined on $E$}, and
$\vm{P}_1^* (E)$ denote the vector polynomial basis on $E$ that is exactly of order $1$. 
For linear elasticity or linear elastodynamics, there are three
rigid-body modes in 2D and six rigid-body modes in 3D. These are also
known as zero-energy modes. Let
\begin{equation}\label{eq:xietazeta}
\xi := \dfrac{x - x_E}{h_E}, \quad \eta := \dfrac{y - y_E}{h_E}, \quad
\zeta := \dfrac{z - z_E}{h_E}.
\end{equation}
For solid continua, we use a suitable scaled vectorial basis set
(first order, $k = 1$) in 2D and 3D. In two dimensions,
\begin{subequations}\label{eq:M}
\begin{align}
\label{eq:M-a}
\vm{M}^\textrm{RBM}(E) &= 
\left[
\begin{Bmatrix} 1 \\ 0 \end{Bmatrix}, \ 
\begin{Bmatrix} 0 \\ 1 \end{Bmatrix}, \ 
\begin{Bmatrix} -\eta \\ \xi \end{Bmatrix} \right], \quad
\widehat{\vm{M}}(E) := \vm{M}_1(E) = \vm{M}^\textrm{RBM}(E) \cup
\left[
\begin{Bmatrix} \eta \\ \xi \end{Bmatrix}, \
\begin{Bmatrix} \xi \\ 0 \end{Bmatrix}, \
\begin{Bmatrix} 0 \\ \eta \end{Bmatrix}
\right] , \\
\intertext{and in three dimensions,}
\label{eq:M-b}
\vm{M}^\textrm{RBM}(E) &= 
\left[ 
\begin{Bmatrix} 1 \\ 0 \\ 0 \end{Bmatrix}, \
\begin{Bmatrix} 0 \\ 1 \\ 0 \end{Bmatrix}, \
\begin{Bmatrix} 0 \\ 0 \\ 1 \end{Bmatrix}, \
\begin{Bmatrix} 0 \\ -\zeta \\ \eta \end{Bmatrix}, \
\begin{Bmatrix} \zeta \\ 0 \\ -\xi \end{Bmatrix}, \
\begin{Bmatrix} -\eta \\ \xi \\ 0 \end{Bmatrix}
\right], \nonumber \\
\widehat{\vm{M}}(E) := 
\vm{M}_1(E) &= \vm{M}^\textrm{RBM}(E) \cup \left[
\begin{Bmatrix} 0 \\ \zeta \\ \eta \end{Bmatrix}, \
\begin{Bmatrix} \zeta \\ 0 \\ \xi \end{Bmatrix}, \
\begin{Bmatrix} \eta \\ \xi \\ 0 \end{Bmatrix}, \
\begin{Bmatrix} \xi \\ 0 \\ 0 \end{Bmatrix}, \
\begin{Bmatrix} 0 \\ \eta \\ 0 \end{Bmatrix}, \
\begin{Bmatrix} 0 \\ 0 \\ \zeta \end{Bmatrix}
\right],
\end{align}        
\end{subequations}
where $\vm{M}^\textrm{RBM}(E)$ contains the rigid body modes (vectors). 
Now, $\widehat{\vm{M}}(E)$ is a basis for $\vm{P}_1(E)$, and we
use $\vm{m}_\alpha$ to refer to a vector element of $\widehat{\vm{M}}(E)$. Further, we use
$\widehat{\vm{M}}^*(E)$ to denote the basis set for a vector field 
on $E$ that contains all
scaled vector monomials that are exactly of order $1$. 

\subsection{Polynomial projection operators}
We use standard Sobolev spaces, $L^2(E)$ and $H^1(E)$, on element $E$
with usual definitions for inner products and norms (seminorms) for
functions that belong to these spaces. Since only $k = 1$ is considered in
this paper, we let $\Pi^0 \equiv \Pi_{E,1}^0$,
$\Pi^\nabla \equiv \Pi_{E,1}^\nabla$,
 and 
$\Pi^\varepsilon \equiv \Pi_{E,1}^\varepsilon$ be the $L^2$, elliptic (Laplacian),
and strain energy (solid continua) projection operators.

\smallskip

\noindent {\bf $\vm{L^2}$ Projection Operator}:
        For a vector field $\vm{v} \in [L^2(E)]^d$, 
        define the $L^2$ orthogonal projection
        operator $\Pi^0 : [L^2(E)]^d \rightarrow 
        \vm{P}_1(E)$ such that
        $\Pi^0 \vm{v}$ is the unique element 
        in $\vm{P}_1(E)$ that satisfies
\begin{equation}\label{eq:L2}
 (\vm{q},\vm{v} - \Pi^0 \vm{v})_E  = 0 
 \quad \forall \vm{q} \in \vm{P}_1(E),
\end{equation}
where $(\cdot,\cdot)_E$ is the standard $L^2$ inner product over $E$.

\smallskip

\noindent {\bf Elliptic (Vector Laplacian) Projection Operator}: 
For a vector-valued function $\vm{v} \in [H^1(E)]^d$ and
the inner product (bilinear form) associated with the vector Laplacian,
$a_E(\vm{u},\vm{v}) = \int_E \nabla \vm{u} : \nabla \vm{v} \, d\vx$,
define the vectorial elliptic orthogonal projection operator 
$\Pi^\nabla : [H^1(E)]^d \rightarrow \vm{P}_1(E)$ such that
$\Pi^\nabla \vm{v}$ is the unique element in $\vm{P}_1(E)$ 
that satisfies
\begin{equation} \label{eq:a-nabla}
\begin{cases}
a_E(\vm{q},\vm{v} - \Pi^\nabla \vm{v}) &= 0 
 \quad \forall \vm{q} \in \vm{P}_1(E), \\
\int_{\partial E}(\vm{v} - \Pi^\nabla \vm{v}) \, dS &= \vm{0} .
\end{cases}
\end{equation}

\smallskip

\noindent {\bf Energy Projection Operator for Solid Continua}: 
Let $\vm{u}, \vm{v} \in [H^1(E)]^3$ be vector-valued functions and
$[\mathbb{P}_{1}(E)]_\textrm{sym}^{3 \times 3}$ be the space of order
$1$ 
polynomial second-order symmetric tensors. 
The internal virtual work for a deformable solid is:
\begin{equation}\label{eq:a-innerproduct}
a_E(\vm{u}, \vm{v}) = \int_E \vm{\sigma}(\vm{u}) : \vm{\varepsilon}  (\vm{v}) \, d \vx,
\end{equation}
where $\vm{u}$ is the trial displacement field, $\vm{v}$ is the test (virtual) displacement field, and
$\vm{\varepsilon}(\vm{v}) $ and $\vm{\sigma}(\vm{u})$ are the virtual strain tensor and Cauchy stress tensor,
respectively.

If $\vm{u}$ is a vector polynomial of order 1, 
then $\vm{\varepsilon} (\vm{u})= \nabla_s (\vm{u}) \in 
[\mathbb{P}_0(E)]_\textrm{sym}^{3 \times 3}$ is a constant symmetric 
second-order tensor ($\nabla_s$ is the symmetric gradient operator)
and so is $\vm{\sigma}(\vm{u})$.   
We define the projection operator
$\Pi^\varepsilon: [H^1(E)]^3 \rightarrow \vm{P}_1(E)$ such 
that $\Pi^\varepsilon \vm{v} \in \vm{P}_1(E)$, whereas
 $\Pi^\varepsilon \vm{\varepsilon} \in [\mathbb{P}_0 (E)]_\textrm{sym}^{3 \times 3}$.
For first order consistency, we require that $(\vm{v} -
    \Pi^\varepsilon \vm{v})$ is 
orthogonal to all vector
polynomials in the energy inner product defined in~\eqref{eq:a-innerproduct}:
\begin{subequations}\label{eq:strain-projection}
\begin{align}
\label{eq:strain-projection-a}
a_E(\vm{m}_\alpha,\vm{v} - \Pi^\varepsilon \vm{v}) &= 0 
 \quad \forall \vm{m}_\alpha \in \widehat{\vm{M}}(E) .
\intertext{Referring to~\eqref{eq:M}, observe 
           that for $\vm{m}_\alpha \in \vm{M}^\textrm{RBM}(E)$, we have
           $ \vm{\varepsilon}(\vm{m}_\alpha) = \vm{0}$, and therefore
           $\vm{\sigma}(\vm{m}_\alpha) = \vm{C} : \vm{\varepsilon}(\vm{m}_\alpha) = \vm{0}$. 
           Hence, for such $\vm{m}_\alpha$, \eqref{eq:strain-projection-a} yields the trivial
           identity $0 = 0$, and hence we provide additional conditions to fully determine the
           coefficients of the vector polynomial that correspond to the rigid-body modes.
           For $k = 1$, we use~\cite{Beirao:2013:VEL}
           }
\label{eq:strain-projection-b}          
\int_{\partial E}  (\vm{v} - \Pi^\varepsilon \vm{v}) \cdot \vm{m}_\alpha
 \, dS &= \vm{0} ,
\end{align}
\end{subequations}
where $\alpha = 1$--$3$ in 2D and $\alpha = 1$--$6$ in 3D.

\section{Virtual Element Method for Solid Continua}\label{sec:vem}
Let $V(E)$ denote the first-order ($k = 1$) virtual element space on 
\suku{a three-dimensional} element
$E$ and $V(\face)$ be the first-order virtual element space
on a particular face $\face$ of $E$.
First, the virtual element space in two dimensions is described and then the virtual element space in three dimensions. Other key elements such as the degrees of freedom, and the construction of the bilinear forms for the internal virtual work and the inertial virtual work are also presented. Finally, the numerical implementation of the method (computation of
consistency and stability parts of the stiffness and mass matrices) is described. The exposition draws on the presentations that appear in the literature.\cite{Beirao:2014:HGV,Cangiani:2015:HSV,Mengolini:2019:EPV,Benvenuti:2019:EVE}

\subsection{Definition of the space
            \texorpdfstring{$\vm{V(F)}$}{TEXT}
            for
            \texorpdfstring{$\vm{d=2}$}{TEXT}}\label{subsec:2D}
Let $F$ be an element in two dimensions. \suku{We point out that in~\sref{subsec:3D},
we use $\face$ (two-dimensional polygon) to represent a face of a three-dimensional 
element $E$.} 
For a vectorial problem, the (enhanced) virtual element space $V(F)$ in two dimensions is defined 
as:~\cite{Ahmad:2013:EPV}
\begin{align}\label{eq:V-2d}
V(F) = \Biggl\{ & \vm{v}^h : \vm{v}^h \in [H^1(F)]^2, \ 
        \Delta \vm{v}^h \in \vm{P}_1(F), \
         \vm{v}^h |_e \in \vm{P}_1(e) \ \forall e \in \partial F,
\nonumber \\
         & \vm{v}^h |_{\partial F} \in [C^0]^2 (\partial F), \, 
         \int_F \vm{v}^h \cdot \vm{q} \, d \vx = 
         \int_F \Pi_F^\nabla \vm{v}^h 
         \cdot \vm{q} \, d \vx \ \ \forall \vm{q} \in
         \vm{P}_{0}^*(F) \cup \vm{P}_{1}^*(F)
          \Biggr\}. 
\end{align}
In two dimensions, a vector-valued function $\vm{v}^h\in V(F)$ is characterized by the following properties:
\begin{enumerate}[(i)]
\item $\vm{v}^h$ is such that each component is continuous and piecewise linear on $\partial F$.
\item $ \Delta \vm{v}^h$ is a linear vector polynomial,
where $ \Delta$ is the vector Laplacian.
\item $\int_F \vm{v}^h \cdot \vm{q} \,d\vx =
                    \int_F \Pi_F^\nabla \vm{v}^h \cdot \vm{q} \, d \vx$ for all linear vector polynomials $\vm{q}$, where
                    $\Pi_F^\nabla$ is the elliptic projector in $F$.                   
\end{enumerate}

Note that the second condition is distinct from the one introduced 
in Beir{\~a}o et al~\cite{Beirao:2013:BPV} for the Poisson problem ($\Delta v^h=0$).
In addition, the third condition states that $\Pi_F^\nabla \vm{v}^h$ 
is also the $L^2$-projection of $\vm{v}^h$ onto linear vector
polynomials for $k = 1$.\cite{Ahmad:2013:EPV}  In essence, the second condition adds three additional degress of freedom in each dimension (Laplacian of a scalar field is an affine function) and the third condition provides an equation to compute them.
These two modifications enable the computations in 3D (integrals over polygonal faces
of the polyhedra)
to be performed using only the values of $\vm{v}^h$ at the vertices of $E$.

\subsection{Definition of the space  
            \texorpdfstring{$\vm{V(E)}$}{TEXT}
            for \texorpdfstring{$\vm{d=3}$}{TEXT}}
            \label{subsec:3D}
For a vectorial problem, the (enhanced) virtual element space $V(E)$ 
in three dimensions is defined as:~\cite{Ahmad:2013:EPV}
\begin{align}\label{eq:V-3d}
V(E) = \Biggl\{ & \vm{v}^h : \vm{v}^h \in [H^1(E)]^3, \, 
         \Delta \vm{v}^h \in \vm{P}_1(E), \,
         \vm{v}^h |_\face \in V(\face) \ \forall \face \in \partial E, \nonumber \\
         &
         \int_E \vm{v}^h \cdot \vm{q} \, d\vx=\int_E \Pi_E^\nabla \vm{v}^h 
         \cdot \vm{q} \,d \vx \ \forall \vm{q} \in
         \vm{P}_{0}^*(E) \cup \vm{P}_{1}^*(E)
          \Biggr\}. 
\end{align}        
In three dimensions, a function $\vm{v}^h \in V(E)$ is characterized by the following properties:
\begin{enumerate}[(i)]
\item $\vm{v}^h$ is continuous and unknown inside $E$, and in
       general also unknown on $\partial E$. 
\item $\vm{v}^h$ is such that its restriction to a face $ \face $ 
       belongs to $ V(\face)$ as defined in~\eqref{eq:V-2d}.
\item  $\Delta \vm{v}^h$ is a linear vector polynomial.
\item  $\int_E \vm{v}^h \cdot \vm{q} \,d\vx =
                    \int_E \Pi_E^\nabla \vm{v}^h \cdot \vm{q} \, d \vx $ for all linear vector polynomials $\vm{q}$.
\end{enumerate}
\suku{
The global virtual element space $V^h \subset [ H_0^1(\Omega)]^3$ is defined using standard assembly procedures as in finite elements:
\begin{equation}\label{eq:Vh}
V^h = \{ \vm{v}^h \in [H_0^1(\Omega)]^3 : \vm{v}^h |_E \in V(E) \textrm{ for all } E \in \Omega^h \}.
\end{equation}
}


\subsection{Degrees of freedom}
To solve a vectorial problem using a first-order virtual element formulation ($k = 1$), 
we take the values of $\vm{v}^h$ at the vertices of the polyhedron as its degrees of freedom (DOFs).
Let us associate a (virtual) canonical
Lagrange basis, $\phi_i(\vx) \  (i = 1,2,\ldots,N_E)$, to each vertex.   The
$\phi_i(\vx)$ possesses the Kronecker-delta property, namely $\phi_i(\vx_j) = \delta_{ij}$. In addition, they satisfy the constant and linear reproducing conditions:
$\sum_{i=1}^{N_E} \phi_i(\vx) = 1$, 
$\sum_{i=1}^{N_E} \phi_i(\vx) \xi_i = \xi$, 
$\sum_{i=1}^{N_E} \phi_i(\vx) \eta_i = \eta$, and
$\sum_{i=1}^{N_E} \phi_i(\vx) \zeta_i = \zeta$. 
The total number of degrees of freedom in each element $E$ is $N_E^{\textrm{dofs}} = 3 N_E$. Let us define the vectorial 
basis function matrix: 
\setcounter{MaxMatrixCols}{20}
\begin{equation} \label{eq:varphi}
\begin{split}
\vm{\varphi } := \vm{\phi} \otimes \vm{I}_3 
= \{ \phi_1  \ \phi_2 \ \dots \ \phi_{N_E} \} \otimes \vm{I}_3
&=
\begin{bmatrix}
\phi_1 & \phi_2 & \dots &  \phi_{N_E}   &   0    &   0   &  \dots  &  0   &   0   &    0   &  \dots  &   0 \\
0 & 0 & \dots &       0   &   \phi_1   &   \phi_2  &  \dots  &  \phi_{N_E}  &   0   &    0   &  \dots  &   0 \\
0 & 0 & \dots &       0   &   0    &   0   &  \dots  &  0   &  \phi_1  & \phi_2   &  \dots  &   \phi_{N_E}
\end{bmatrix} \\
&:= 
\begin{bmatrix}
\vm{\varphi}_1 & \vm{\varphi}_2 & \dots & \vm{\varphi}_{N_E} & 
\vm{\varphi}_{N_E+1} & \vm{\varphi}_{N_E+2} & \dots & \vm{\varphi}_{2N_E} & 
\vm{\varphi}_{2N_E+1} & \vm{\varphi}_{2N_E+2} & \dots & \vm{\varphi}_{3N_E} 
\end{bmatrix},
\end{split}
\end{equation}
where $ \vm{\varphi}_{i} = \{ \phi_i \ 0 \ 0 \}^T$, 
$\vm{\varphi}_{N_E+i} = \{ 0 \ \phi_i \ 0 \}^T$, and
$\vm{\varphi}_{2N_E+i} = \{ 0 \ 0 \ \phi_i\}^T$ for $i = 1,\dots,N_E$, are the $3N_E$
vectorial basis functions. This choice can be seen as a natural generalization of the scalar problem to the vectorial case.
We can write the interpolant of a vector-valued function $\vm{v}(\vx)$ in $E$ as:
\begin{equation}\label{eq:vh1} 
\vm{v}^h(\vx) = \sum_{i=1}^{3 N_E} \vm{\varphi}_i (\vx) v_i  ,
\end{equation}
where each vectorial basis function is multiplied by a scalar coefficient. We can 
express~\eqref{eq:vh1} as
\begin{equation}\label{eq:vh2}
\vm{v}^h(\vx) = \sum_{i=1}^{3 N_E} \vm{\varphi}_i(\vx) \, \texttt{dof}_i (\vm{v^h}),
\end{equation}
where the operator $\texttt{dof}_i(\vm{v}^h)$ extracts the $i$-th degree of freedom (scalar) of the vector field $\vm{v}^h$.

\subsection{Discrete generalized eigenproblem}\label{subsec:discrete-eig}
From the weak form of the continuous problem in~\eqref{eq:weak},
we can write the weak form for the VEM as:~find
$\vm{u}^h \in V^h \subset V $ and $\lambda \in \Re_+$ such that
\begin{align}\label{eq:weak-discrete}
  a^h(\vm{u}^h,\vm{v}^h) = \lambda \, b(\vm{u}^h,\vm{v}^h) \quad 
  \forall \vm{v}^h \in V^h \subset V,
\end{align}
where $a^h(\cdot,\cdot)$ and $b^h(\cdot,\cdot)$ are the 
virtual element bilinear forms that approximate the exact bilinear
forms $a(\cdot,\cdot)$ and $b(\cdot,\cdot)$.

We define the discrete bilinear forms $a^h(\vm{u}^h,\vm{v}^h)$ 
and $b^h(\vm{u}^h,\vm{v}^h)$ for
$\vm{u}^h, \vm{v}^h \in V^h$ as the sum of elemental contributions
\begin{align}
 a^h(\vm{u}^h,\vm{v}^h)=\sum_{E \in {\cal T}^h } 
 a_E^h (\vm{u}^h,\vm{v}^h), \quad
 b^h(\vm{u}^h,\vm{v}^h)=\sum_{E \in {\cal T}^h } 
 b_E^h (\vm{u}^h,\vm{v}^h),
\end{align}
where the weak form on element $E$ reads: 
find $\vm{u}^h \in V(E)$ and $\lambda_E \in \Re_+$ such that
\begin{subequations}\label{eq:weak-E}
\begin{align}
\label{eq:weak-E-a}
a_E^h(\vm{u}^h,\vm{v}^h) &= \lambda_E \, b_E^h(\vm{u}^h,\vm{v}^h) \quad 
  \forall \vm{v}^h \in V(E), \\
\label{eq:weak-E-b}  
a_E^h(\vm{u}^h,\vm{v}^h ) &= \int_E \vm{\sigma}(\vm{u}^h) :
\vm{\varepsilon}(\vm{v}^h) \, d\vx, \quad
b_E^h(\vm{u}^h,\vm{v}^h) = \int_E \vm{u}^h \cdot \vm{v}^h \, d\vx .
\end{align}
\end{subequations}

Recall that $\{\vm{\varphi}_i\}_{i=1}^{3N_E}$ are the 
canonical vectorial basis function on $E$. 
On expanding the trial and test functions as a linear combination of
these basis functions and substituting them in~\eqref{eq:weak-E},
and using the arbitrariness of the test coefficients, we obtain the
following generalized eigenproblem:
\begin{subequations}\label{eq:Kd}
  \begin{align}
  \vm{K} \vm{d} &= \lambda \, \vm{M} \vm{d}, \quad
  \vm{K} = \sum_{E \in {\cal T}^h} \vm{K}_E, \
  \vm{M} = \sum_{E \in {\cal T}^h} \vm{M}_E, \\
  \vm{K}_E &= a_E^h(\vm{\varphi},\vm{\varphi}), \quad
  \vm{M}_E = b_E^h(\vm{\varphi},\vm{\varphi}) ,
  \end{align}
\end{subequations}
where $\vm{K}_E$ and $\vm{K}$ are the element and global stiffness matrices, 
$\vm{M}_E$ and $\vm{M}$ are the element and global mass matrices,
and the pair $\{\lambda,\vm{d}\}$ is the eigensolution. 
\suku{In the VEM, $a_E^h(\cdot,\cdot)$ and
$b_E^h(\cdot,\cdot)$ consist of a consistency term and a stability term. These
terms for the stiffness matrix are presented in~\eqref{eq:K_E} and those for
the mass matrix appear in~\eqref{eq:M_E}.}
Note 
that~\eqref{eq:weak-E} is solved to compute the maximum element
eigenvalue, which is used to bound the critical time step as shown in~\sref{subsec:timestep}.

For convergence, the bilinear forms $a_E^h(\cdot,\cdot)$ and $b_E^h(\cdot,\cdot)$
must satisfy the following properties:\cite{Beirao:2013:BPV,Vacca:2015:VEM}
\begin{description}
\item[$\bullet$] {\emph{linear consistency}}: for all $\vm{v}^h \in V(E)$ and linear vector polynomials
 $\vm{q} \in \vm{P}_1(E)$  it holds that
   \begin{subequations}\label{eq:linear-consistency}
  \begin{align} 
    \label{eq:linear-consistency-a}
    a_E^h(\vm{q},\vm{v}^h) &= a_E(\vm{q},\vm{v}^h), \\
    \label{eq:linear-consistency-b}
    b_E^h(\vm{q},\vm{v}^h) &= (\vm{q},\vm{v}^h)_E .
  \end{align}
  \end{subequations}
\item[$\bullet$] {\emph{stability}}: there exists four positive constants
  $\alpha_*$,$\, \alpha^*$,$\, \beta_*$,$\, \beta^*$, independent of $h$ and $E$, such that
  \begin{subequations}\label{eq:stability}
  \begin{align}
    \label{eq:stability-a}
     \alpha_* a_E(\vm{v}^h,\vm{v}^h) \le a_E^h(\vm{v}^h,\vm{v}^h) 
     & \le \alpha^* a_E(\vm{v}^h,\vm{v}^h) \quad \forall \vm{v}^h \in V(E), \\
     \label{eq:stability-b}
     \beta_* (\vm{v}^h,\vm{v}^h)_E \le b_E^h(\vm{v}_h,\vm{v}_h) 
     & \le \beta^* (\vm{v}^h,\vm{v}^h)_E \quad \forall \vm{v}^h \in V(E) .
     \end{align}
     \end{subequations}
\end{description}

\section{Numerical Implementation}\label{sec:implementation}
We derive the 
discrete system to compute the matrix representation of the energy projector
$\Pi^\epsilon$ and the $L^2$ projector
$\Pi^0$.  Expressions are presented for the 
element stiffness matrix in three dimensions, and for the consistent and lumped mass matrices in both 2D and 3D.
Finally, the main elements of the homogeneous numerical integration method are 
discussed.

\subsection{Computation of energy projection matrices}
Referring to~\eqref{eq:varphi}, let
\begin{equation}\label{eq:Pivarphi_i}
    \Pi^\varepsilon \vm{\varphi}_i = 
    \sum_{\beta = 1}^{12} \vm{m}_\beta  s_\beta^i 
    = \widehat{\vm{M}} \vm{s}^i \ \ 
    (i = 1,2,\dots,3N_E)
\end{equation}
be the projection of the vectorial canonical basis functions onto the scaled vector polynomial basis functions. In~\eqref{eq:Pivarphi_i}, $s_\beta^i$ are unknown
coefficients, 
 $\widehat{\vm{M}} := 
\widehat{\vm{M}}(E)$ is given in~\eqref{eq:M},
and we have suppressed the spatial dependence of
$\vm{\varphi}_i$ and $\vm{m}_\alpha$. Now, we can write 
\begin{equation}\label{eq:Pivarphi}
    \Pi^\varepsilon \vm{\varphi} =
    \Pi^\varepsilon \bigl\{ \vm{\varphi}_1 \  \vm{\varphi}_2  \dots
                       \vm{\varphi}_{3N_E} \bigr\} =
                       \sum_{\beta = 1}^{12} \vm{m}_\beta
                       \bigl\{ s_\beta^1 \ s_\beta^2 
                           \dots s_\beta^{3N_E} \bigr\}
                           := \widehat{\vm{M}} \vm{S} 
                           := \widehat{\vm{M}} \vm{\Pi}_*^\varepsilon ,
\end{equation}
where $\vm{\Pi}_*^\varepsilon := \vm{S}$ is the 
matrix representation of the \suku{projection of the canonical shape functions}
with respect to the scaled vector monomial basis set $\widehat{\vm{M}}$.

Define the dof-matrix, $\vm{D}$, with entries
$\vm{D}_{i\alpha} =  \texttt{dof}_i(\vm{m}_\alpha)$ as:
\begin{equation}\label{eq:D}
\vm{D} = 
    \begin{bmatrix}
    \texttt{dof}_1 (\vm{m}_1) & \texttt{dof}_1 (\vm{m}_2) &  \dots & \texttt{dof}_1 (\vm{m}_{12}) \\ 
    \texttt{dof}_2 (\vm{m}_1) & \texttt{dof}_2 (\vm{m}_2) &  \dots & \texttt{dof}_2 (\vm{m}_{12}) \\ 
    \dots & \dots & \dots & \dots \\
    \texttt{dof}_{3N_E} (\vm{m}_1) & \texttt{dof}_{3N_E} (\vm{m}_2) &  \dots & \texttt{dof}_{3N_E} (\vm{m}_{12}) \\ 
    \end{bmatrix} ,
\end{equation}
where the operator $\texttt{dof}_i(\cdot)$ is defined in~\eqref{eq:vh2}. 
On using~\eqref{eq:M}, \eqref{eq:varphi} and~\eqref{eq:D}, the constant and linear reproducing conditions can be written as: 
\begin{equation}\label{eq:linear_consistency}
    \widehat{\vm{M}} = \vm{\varphi} \vm{D} = 
    \vm{I} \otimes \{ \phi_1 \ \phi_2 \dots \ \phi_{N_E} \}
    \begin{bmatrix}
    \vm{1} & \vm{0} & \vm{0}  & \vm{0} & \zetavec &
    -\etavec & \vm{0} & \zetavec & \etavec &
    \xivec & \vm{0} & \vm{0} \\
    \vm{0} & \vm{1} & \vm{0}  & - \zetavec &
    \vm{0} & \xivec & \zetavec & \vm{0} &
     \xivec & \vm{0} & \etavec & \vm{0} \\
    \vm{0} & \vm{0} & \vm{1}  & \etavec &
    - \xivec & \vm{0} & \etavec & \xivec &
    \vm{0} & \vm{0} & \vm{0} & \zetavec
    \end{bmatrix} ,
\end{equation}
where $\vm{1}$ is a $N_E \times 1$ column vector of ones, $\vm{0}$ is a $N_E \times 1$ column vector of zeros, and
$\xivec = \{ \xi_1 \ \xi_2 \ \dots \ \xi_{N_E} \}^T$,
$\etavec = \{ \eta_1 \ \eta_2 \ \dots \ \eta_{N_E} \}^T$,
$\zetavec = \{ \zeta_1 \ \zeta_2 \ \dots \ \zeta_{N_E} \}^T$
are the nodal values of the affine functions $\xi$, $\eta$ and 
$\zeta$, respectively, which are defined in~\eqref{eq:xietazeta}.

On using~\eqref{eq:linear_consistency}, we have
\begin{subequations}
\begin{align}
    \Pi^\varepsilon \vm{\varphi} &= \widehat{\vm{M}} \vm{\Pi}_*^\varepsilon = \vm{\varphi} \vm{D} 
    \vm{\Pi}_*^\varepsilon ,
\intertext{and therefore}
\vm{\Pi}^\varepsilon &:= \vm{D} \vm{\Pi}_*^\varepsilon
\end{align}
\end{subequations}
is the matrix representation of the projection in the
$\vm{\varphi}$ basis set.

\suku{We define a mean discrete $L^2$ inner product on the element $E$} as:
\begin{equation}\label{eq:P_0}
    P_0(\vm{u}^h,\vm{v}^h) :=
    \frac{1}{N_E} \sum_{k=1}^{N_E} 
    \vm{u}^h (\vx_k) \cdot  \vm{v}^h (\vx_k) 
    = 
    \frac{1}{N_E} \sum_{i=1}^{3N_E} 
    \texttt{dof}_i (\vm{u}^h) \, \texttt{dof}_i ( \vm{v}^h ) .
\end{equation}
Now, from~\eqref{eq:strain-projection}, the variational problem to determine the projector can be written as
\begin{subequations}\label{eq:vem-projection}
\begin{align}
\label{eq:vem-projection-a}
a_E^h(\vm{m}_\alpha, \Pi^\varepsilon \vm{v}^h ) &= 
a_E^h(\vm{m}_\alpha, \vm{v}^h ) 
 \quad \textrm{for } \vm{m}_\alpha \in \widehat{\vm{M}} \
 (\alpha = 7,8,\dots,12), \\
\label{eq:vem-projection-b}  
P_0(\vm{m}_\alpha, \Pi^\varepsilon \vm{v}^h)_{\partial E} 
&= P_0(\vm{m}_\alpha, \vm{v}^h)_{\partial E} \quad
\textrm{for } \vm{m}_\alpha \in \widehat{\vm{M}} \
 (\alpha = 1,2,\dots,6) .
\end{align}
\end{subequations}
On substituting $\vm{v}^h = \vm{\varphi}_i$ ($i = 1,2,\ldots 3N_E$)
in~\eqref{eq:vem-projection} and
using~\eqref{eq:P_0}, we have
\begin{subequations}\label{eq:vem-proj}
\begin{align}
\label{eq:vem-proj-a}
a_E^h(\vm{m}_\alpha, \Pi^\varepsilon \vm{\varphi}_i ) &= 
a_E^h(\vm{m}_\alpha, \vm{\varphi}_i ) 
 \quad \textrm{for } \vm{m}_\alpha \in \widehat{\vm{M}} \
 (\alpha = 7,8,\dots,12), \\
\label{eq:vem-proj-b}  
\frac{1}{N_E} \sum_{j=1}^{3N_E} 
    \texttt{dof}_j (\vm{m}_\alpha) \, \texttt{dof}_j ( \Pi^\epsilon \vm{\varphi}_i ) 
    &= 
    \frac{1}{N_E} \sum_{j=1}^{3N_E} 
    \texttt{dof}_j (\vm{m}_\alpha) \, \texttt{dof}_j ( \vm{\varphi}_i ) 
 \quad
\textrm{for } \vm{m}_\alpha \in \widehat{\vm{M}}  \
 (\alpha = 1,2,\dots,6) .
\end{align}
\end{subequations}
Equation~\eqref{eq:vem-proj-b} equates the average value of
$\vm{\varphi}_i$ and that of its projection along the coordinate
directions (rigid-body translations for $\alpha = 1,2,3$) and
likewise for the infinitesimal rigid-body rotations about the $x$-, 
$y$-, and $z$-axes $(\alpha = 4,5,6)$.

On using~\eqref{eq:Pivarphi_i}, \eqref{eq:vem-proj-a} becomes
\begin{align*}
\sum_{\beta = 1}^{12} a_E^h ( \vm{m}_\alpha, \vm{m}_\beta )
s_\beta^i 
&= a_E^h(\vm{m}_\alpha.\vm{\varphi}_i) \quad (\alpha = 7,8, \dots, 12) \\
&= \int_E \vm{\sigma} (\vm{m}_\alpha) :  
    \vm{\varepsilon} (\vm{\varphi}_i)  \, d\vx \\ 
&=  \int_E    \nabla \cdot \bigl( \vm{\varphi}_i \cdot
    \vm{\sigma}(\vm{m}_\alpha) \bigr) \, d\vx
     - \int_E \bigl( 
     \cancel{ \nabla \cdot \vm{\sigma} ( \vm{m}_\alpha ) }
     \bigr) \cdot \vm{\varphi}_i \, d\vx \\
     &= \int_{\partial E}  \vm{t}(\vm{m}_\alpha) \cdot \vm{\varphi}_i
     \, dS 
     =      \sum_{\face \subset \partial E}
     \vm{t}_\face  (\vm{m}_\alpha) \cdot \int_\face \Pi_\face^\varepsilon 
     \vm{\varphi}_i \, dS ,
\end{align*}
where the divergence theorem has been invoked and we have used the fact that $\nabla \cdot \vm{\sigma}(\vm{m}_\alpha) = \vm{0}$ since
$\vm{m}_\alpha$ are affine fields. In addition,
$\vm{t}_\face(\vm{m}_\alpha)$ is the traction vector on $\face$ and $\Pi_\face^\varepsilon$ is the projection operator on $\face$.
The three vectorial basis functions that are associated with node $K$ can be expressed as $\vm{\varphi}_K := 
\phi_K \vm{I}_3$ ($\phi_K$ is the scalar canonical basis function). 
Since finite elements with only
triangular faces are considered in~\sref{sec:results}, then
$\Pi_\face^\varepsilon \vm{\varphi}_K = \vm{\varphi}_K = 
\phi_K \vm{I}_3$ ($\phi_K$ is an
affine function on a triangle, \suku{$\triangle$}) and a one-point
Gauss rule delivers exact integration on each face. If so, then
\begin{equation}\label{eq:Pi-2}
    \sum_{\beta = 1}^{12} a_E^h ( \vm{m}_\alpha, \vm{m}_\beta )
\vm{s}_\beta^K = \sum_{\face = \marksymbol{triangle}{black} } \vm{t}_\face  (\vm{m}_\alpha)  
                 \cdot \int_{\face}  \phi_K \vm{I}_3 \, dS 
               = \frac{1}{3} \sum_{\face = \marksymbol{triangle}{black} } \vm{t}_\face  (\vm{m}_\alpha) |\face|  \quad (\alpha = 7,8,\dots,12).
\end{equation}
On substituting for $\Pi^\varepsilon \vm{\varphi}_i$ from~\eqref{eq:Pivarphi_i} in~\eqref{eq:vem-proj-b}, we obtain
\begin{equation}\label{eq:Pi-1}
    \frac{1}{N_E} \sum_{j=1}^{3N_E} \sum_{\beta=1}^{12}
    \texttt{dof}_j (\vm{m}_\alpha) \, \texttt{dof}_j 
    ( \vm{m}_\beta ) \, \vm{s}_\beta^K 
    = 
    \frac{1}{N_E} \sum_{j=1}^{3N_E} 
    \texttt{dof}_j (\vm{m}_\alpha) \, \texttt{dof}_j ( \vm{\varphi}_K ) 
 \quad
 (\alpha = 1,2,\dots,6).
\end{equation}
Using Voigt notation, let $\widetilde{\vm{B}} = \vm{\varepsilon}(\widehat{\vm{M}})$ be the $6 \times 12$
matrix that represents the symmetric gradient of the column vectors
of $\widehat{\vm{M}}$ in~\eqref{eq:M} and $\vm{C}$ be the $6 \times 6$
constitutive matrix for an isotropic linear elastic solid. Also,
let $\vm{\sigma}^* = \vm{C} \widetilde{\vm{B}}$ be the $6 \times 12$
matrix representation of the stress tensor of the
displacement fields in $\widehat{\vm{M}}$; then the traction on a face $\face$ is 
$\vm{t}_\face = \vm{n}_\face \cdot \vm{\sigma}^*.$
In~\eqref{eq:Pi-2}, $a_E^h(\vm{m}_\alpha,\vm{m}_\beta) = 
\int_E \vm{\sigma} (\vm{m}_\alpha) : \vm{\varepsilon} (\vm{m}_\beta) 
\, d\vx = \vm{\sigma} ( \vm{m}_\alpha) : \vm{\varepsilon} (\vm{m}_\beta)
|E| = \widetilde{\vm{B}}^T \vm{C} \widetilde{\vm{B}} |E|$.
Further, on noting that $\texttt{dof}_i (\vm{m}_\alpha) = \vm{D}_{i\alpha}$,
$\texttt{dof}_j ( \vm{\varphi}_i ) = \delta_{ij}$ and recalling
that $\vm{S} := \vm{\Pi}_*$ is the matrix of unknown coefficients,
the linear system that stems from~\eqref{eq:Pi-1} and~\eqref{eq:Pi-2} can 
be written as:
\begin{subequations}
\begin{align}
    \vm{G} \vm{\Pi}_* &= \widehat{\vm{B}},  \quad
    \vm{\Pi}_* = \vm{G}^{-1} \widehat{\vm{B}}, 
    \intertext{where}
    \widetilde{\vm{G}} = \widetilde{\vm{B}}^T \vm{C} \widetilde{\vm{B}} 
    |E| ,  \quad
            \vm{G} &= \widetilde{\vm{G}}, \quad
            \vm{G}(1\!:\!6,:) =
            \frac{ [\vm{D}^T \vm{D}] (1\!:\!6,:) }{N_E}, \\
            \widehat{\vm{B}}(1\!:\!6,:) &= 
            \frac{[ \vm{D}(:,1\!:\!6) ]^T }{N_E}, \quad 
            \widehat{\vm{B}}(7\!:\!12,:) = 
            \frac{1}{3} \sum_{\face = \marksymbol{triangle}{black} \subset \partial E}  
            [ (\vm{\sigma}^*)^T \vm{n}_f ] (7\!:\!12,:) \, |\face|.
\end{align}
\end{subequations}

\subsection{Element stiffness matrix}
The expression for the element stiffness matrix is:
\begin{align*}
\vm{K}_E &= a_E^h(\vm{\varphi}, \vm{\varphi}) \\
         &= a_E^h \bigl(\Pi^\varepsilon \vm{\varphi} + (1-\Pi^\varepsilon)
         \vm{\varphi}, \Pi^\varepsilon \vm{\varphi} + 
         (1-\Pi^\varepsilon) \vm{\varphi} \bigr) \\
        &= a_E^h(\Pi^\varepsilon \vm{\varphi}, 
                  \Pi^\varepsilon \vm{\varphi}) 
           + a_E^h\bigl( (1 - \Pi^\varepsilon) \vm{\varphi}, 
                   (1 - \Pi^\varepsilon) \vm{\varphi} \bigr)    \\
        &= a_E^h(\widehat{\vm{M}}\vm{\Pi}_*^\varepsilon,     
                 \widehat{\vm{M}}\vm{\Pi}_*^\varepsilon )   
            + \suku{ a_E^h\bigl( \vm{\varphi} ( \vm{I} - \vm{\Pi}^\varepsilon) , 
                    \vm{\varphi} (\vm{I} - \vm{\Pi}^\varepsilon) \bigr) }    \\          
        &= ( \vm{\Pi}_*^\varepsilon)^T     
           a_E^h(\widehat{\vm{M}} ,  \widehat{\vm{M}} ) 
           \, \vm{\Pi}_*^\varepsilon                   
           + (\vm{I} - \vm{\Pi}^\varepsilon)^T \,  a_E^h(\vm{\varphi},
              \vm{\varphi}) \
             (\vm{I} - \vm{\Pi}^\varepsilon ) ,   
\end{align*}
where we recall that $\vm{\Pi}^\varepsilon = \vm{D} \vm{\Pi}_*^\varepsilon$,
 $\vm{I}$ is the $\suku{dN_E \times dN_E}$ identity matrix,
and the third equality is reached since the cross terms vanish due to the orthogonality condition in~\eqref{eq:strain-projection-a}. The first (consistency) term in $\vm{K}_E$ is computable whereas 
$a_E^h(\vm{\varphi},\vm{\varphi})$ in the second (stability) term is not computable. To ensure coercivity of $\vm{K}_E$, we approximate this term by a diagonal matrix, $\vm{S}_E^d$,
with the appropriate magnitude (scaling). Effective choices for this diagonal matrix (so-called D-recipe)
have been investigated,~\cite{Beirao:2017:HOV,Dassi:2018:EHO} and tailored for
elastodynamics in Park et al.\cite{Park:2019:ONM} 
The entries in the diagonal matrix are chosen as:
$(S_E^d)_{ii} = 
\max \bigl( \alpha_0 \, \texttt{tr} \, (\vm{C})/m, \, (K_E^c)_{ii} 
\bigr)$, 
where
$m = 3$ in 2D and $m = 6$ in 3D, and $\alpha_0 = 1$
in 2D and $\alpha_0 = h_E$ is used in 3D.
On setting $a_E^h(\vm{\varphi},\vm{\varphi}) \approx
\vm{S}_E^d$, both the consistency and stability matrices are now computable:
\begin{equation}
    \vm{K}_E^c = ( \vm{\Pi}_*^\varepsilon )^T 
                    a_E^h( \widehat{\vm{M}} , \widehat{\vm{M}} ) \
                    \vm{\Pi}_*^\varepsilon
                  = ( \vm{\Pi}_*^\varepsilon )^T 
                    \left(
                    \int_E \vm{\sigma}(\widehat{\vm{M}} ) :
                    \vm{\varepsilon}(\widehat{\vm{M}} ) \, d\vx
                    \right) 
                     \vm{\Pi}_*^\varepsilon
                   = ( \vm{\Pi}_*^\varepsilon )^T \, \widetilde{\vm{G}} \,
                            \vm{\Pi}_*^\varepsilon, \quad
    \vm{K}_E^s    = ( \vm{I} - \vm{\Pi}^\varepsilon )^T \, \vm{S}_E^d 
                     ( \vm{I} - \vm{\Pi}^\varepsilon ) .
\end{equation}
Hence, the element stiffness matrix is:
\begin{equation}\label{eq:K_E}
    \vm{K}_E = \vm{K}_E^c + \vm{K}_E^s  = 
    ( \vm{\Pi}_*^\varepsilon )^T \, \widetilde{\vm{G}} \,
                            \vm{\Pi}_*^\varepsilon + ( \vm{I} - \vm{\Pi}^\varepsilon )^T \, \vm{S}_E^d 
                     ( \vm{I} - \vm{\Pi}^\varepsilon ) .
\end{equation}

\subsection{Element mass matrix}
We now use the $L^2$ projection operator to derive the element mass matrix. 
For $k = 1$, the $L^2$ and the elliptic projector coincide:
$\Pi^0 = \Pi^\nabla$.\cite{Ahmad:2013:EPV}  Both satisfy the equality
conditions given in~\eqref{eq:V-3d}, and hence the equivalence. So on
computing the elliptic projection operator we can use its matrix representation to construct the element mass
matrix. The $\Pi^\nabla$ elliptic projector has appeared extensively in the VEM literature; for the sake of completeness, we present the derivation with the most important steps.

Again, let $\vm{\varphi}_K = \phi_K \vm{I}_d$ represent $d$ vectorial basis functions that are associated with node $K$.
Referring to~\eqref{eq:varphi}, 
let the elliptic (vector Laplacian) projection of the vectorial basis functions be represented as: 
\begin{equation}\label{eq:Pivarphi-M}
    \Pi^\nabla \vm{\varphi}_K = \Pi^\nabla ( \phi_K \otimes \vm{I}_d )
    = \Pi^\nabla \phi_K \otimes \vm{I}_d
    = \left( \sum_{\beta = 1}^{d+1} m_\beta \, s_\beta^K \right)
       \otimes \vm{I}_d, 
    \ \ (K = 1,2,\dots,N_E), \ \  m_\beta \in \widehat{\vm{m}}_0
    ,
\end{equation}
where $\widehat{\vm{m}}_0 = \{1, \ \xi, \ \eta \}$ if
$d = 2$ and $\widehat{\vm{m}}_0 = \{1, \ \xi, \ \eta, \ \zeta\}$ 
if $d = 3$, 
and $\vm{S}_0 := \{ \vm{s}^1 \ \vm{s}^2 \ \dots \ \vm{s}^{N_E} \}$ is
the matrix of unknown coefficients. Define the dof-matrix,
$\vm{D}_0 = \texttt{dof}_i (m_\alpha)$, as
\begin{equation}
    \vm{D}_0 = \bigl[
     \vm{1} \  \xivec  \ \etavec 
     \bigr]
    \ \ (d = 2), \quad
    \vm{D}_0 = \bigl[
     \vm{1} \  \xivec  \ \etavec \ \zetavec 
     \bigr] \ \ (d = 3).
\end{equation}
On using $\vm{q} = m_\alpha \otimes \vm{I}_d$ and $\vm{v} :=
\vm{v}^h = \vm{\varphi}_K 
= \phi_K \otimes \vm{I}_d$ in~\eqref{eq:a-nabla}, 
the variational problem 
to determine the elliptic projector becomes
\begin{subequations}\label{eq:vem-elliptic-projection}
\begin{align}
\label{eq:vem-elliptic-projection-a}
a_E^h\bigl(\vm{I}_d \otimes {m}_\alpha, \Pi^\nabla (\phi_K \otimes
\vm{I}_d) \bigr) &= 
a_E^h(\vm{I}_d \otimes {m}_\alpha, \phi_K \otimes \vm{I}_d ) ,
 \quad \textrm{for } m_\alpha \in \widehat{\vm{m}_0}\backslash 1 \\
 \label{eq:vem-elliptic-projection-b}
P_0\bigl(\vm{I}_d \otimes 
{m}_\alpha , \Pi^\nabla (\phi_K \otimes \vm{I}_d) \bigr)_{\partial E} 
&= P_0 (\vm{I}_d \otimes m_\alpha, 
\phi_K \otimes \vm{I}_d)_{\partial E} \quad
\textrm{for } m_\alpha = 1 ,
\end{align}
\end{subequations}
where $K = 1,2,\dots,N_E$,
and the projector $P_0(\cdot,\cdot)$ is defined in~\eqref{eq:P_0}.
On using~\eqref{eq:Pivarphi-M} and~\eqref{eq:P_0}, we can
rewrite~\eqref{eq:vem-elliptic-projection} as
\begin{subequations}\label{eq:scalar-laplacian}
\begin{align}
\vm{I}_d \otimes 
a_E^h \left( m_\alpha, \sum_{\beta} m_\beta s_\beta^K  \right)  \otimes \vm{I}_d &= \vm{I}_d \otimes
a_E^h(m_\alpha, \phi_K ) \otimes \vm{I}_d ,
 \quad \textrm{for }  m_\alpha \in \widehat{\vm{m}_0}\backslash 1 \\
\vm{I}_d \otimes \left(
\frac{1}{N_E} \sum_{j=1}^{N_E} \sum_{\beta}
    \texttt{dof}_j (m_\alpha) \, \texttt{dof}_j ( m_\beta )
    \, \vm{s}_\beta^K \right) \otimes \vm{I}_d
    &= \vm{I}_d \otimes \left(
    \frac{1}{N_E} \sum_{j=1}^{ N_E} \texttt{dof}_j (m_\alpha) \,
     \texttt{dof}_j ( \phi_K ) \right) \otimes \vm{I}_d
     \quad
\textrm{for } m_\alpha = 1 .
\end{align}
\end{subequations}
In~\eqref{eq:scalar-laplacian}, 
the bilinear form, $a_E(u,v) = \int_E \nabla u \cdot
\nabla v \, d\vx$, is
associated with the scalar
Laplace operator. From~\eqref{eq:scalar-laplacian}, we now obtain
\begin{subequations}\label{eq:Pi-nabla}
\begin{align}
\sum_{\beta} \int_E \nabla m_\alpha \cdot \nabla m_\beta \, d\vx \, s_\beta^K 
&= \int_E \nabla m_\alpha \cdot \nabla \phi_K \, d\vx  \nonumber \\
&= \int_E \nabla \cdot (\phi_K \nabla m_\alpha) \, d\vx -
  \int_E \phi_K 
   \cancel{\nabla^2 m_\alpha}  \, d\vx \nonumber \\
&= \int_{\partial E} \phi_K ( \nabla m_\alpha \cdot \vm{n}) \, dS 
   \nonumber \\
 \label{eq:Pi-nabla-a}  
&= \sum_{\face \subset \partial \Omega} 
  \int_{\face} \Pi_\face^\nabla \phi_K ( \nabla m_\alpha \cdot \vm{n}_\face) \, dS 
   \quad \textrm{for }  m_\alpha \in \widehat{\vm{m}_0}\backslash 1, \\
   \intertext{and}
   \label{eq:Pi-nabla-b}
\frac{1}{N_E} \sum_{j=1}^{N_E} \sum_{\beta}
    \texttt{dof}_j  ( m_\alpha )
     \texttt{dof}_j  ( m_\beta ) \, \vm{s}_\beta^K 
    &= 
    \frac{1}{N_E} \sum_{j=1}^{N_E} 
     \delta_{jK} = \frac{1}{N_E} \quad
\textrm{for } m_\alpha = 1 ,
\end{align}
\end{subequations}
where $\nabla^2 m_\alpha = 0$ is used to reach~\eqref{eq:Pi-nabla-a},
$\texttt{dof}_j(m_\alpha) = 1$ for $m_\alpha = 1$,
and
$\texttt{dof}_j(\phi_K) = \delta_{jK}$. Now, 
\begin{equation}
\nabla \widehat{\vm{m}}_0 = \begin{bmatrix}
 0 & \vm{0}^T \\
 \vm{0} & \dfrac{\vm{I}_d}{h_E} 
\end{bmatrix} ,
\end{equation}
where
$\vm{0}$ is a $d$-dimensional column vector of zeros. Let
$\vm{S}_0 = \{ \vm{s}^1 \ \vm{s}^2 \ \dots \ \vm{s}^{N_E} \} $
be the matrix of unknown coefficients.  Then, using the expressions
for $\vm{D}_0$ and $\nabla \widehat{\vm{m}}_0$ 
in~\eqref{eq:Pi-nabla}, we solve for the projection matrix 
$\vm{S}_0$ using the
following linear system of equations:
\begin{subequations}\label{eq:GS_0}
\begin{align}
\vm{G}_0 \vm{S}_0 &= \widehat{B}_0, \quad
\vm{S}_0 = \vm{G}_0^{-1} \widehat{B}_0, \\
\vm{G}_0(1.:) &= \dfrac{\vm{D}_0^T \vm{D} (1,:) }{N_E}, \quad
\vm{G}_0(2\!:\!d\!+\!1,:) = \dfrac{|E|}{h_E^2} \vm{I}_d, \\
\widehat{\vm{B}}_0(1,:) &= \dfrac{\vm{D}_0^T(1,:) }{N_E} = 
\dfrac{\vm{1}}{N_E}, \quad
\widehat{\vm{B}}_0(2\!:\!d\!+\!1,:) = \sum_{\face \subset \partial E} 
\frac{\vm{n}_\face |\face|}{d h_E} ,
\end{align}
\end{subequations}
where $\vm{1}$ is a $N_E$-dimensional row vector of ones
and we have assumed that the boundary faces $\face$ are simplices
in $\Re^d$ so that $\int_\face \Pi_\face^\nabla \phi_K \, dS 
= \int_\face \phi_K \, dS
= |\face|/d$. 
Now, we can write the matrix representations of the $L^2$ projection
operators
(identical to the $\Pi^\nabla$ projection operators) with respect to the scaled vectorial basis,
$\widehat{\vm{m}}_0 \otimes \vm{I}_d$, and the
canonical vectorial basis, $\vm{\phi} \otimes \vm{I}_d$, as
\begin{equation}
    \vm{\Pi}_*^0 = \vm{S}_0 \otimes \vm{I}_d \ \ \textrm{and } \
    \vm{\Pi}^0   =  [ \vm{D}_0 \otimes \vm{I}_d ] \, \vm{\Pi}_*^0,
\end{equation}
respectively.  Proceeding as we did for the stiffness matrix, on
using~\eqref{eq:weak-b} we can write the element mass matrix as
\begin{align}\label{eq:M_E}
    \vm{M}_E &= b_E^h(\vm{\varphi}, \vm{\varphi})  \nonumber \\
         &= b_E^h \bigl(\Pi^0 \vm{\varphi} + (1-\Pi^0)
            \vm{\varphi}, \Pi^0 \vm{\varphi} + 
            (1 - \Pi^0) \vm{\varphi} \bigr) \nonumber \\
        &= b_E^h(\Pi^0 \vm{\varphi}, 
                 \Pi^0 \vm{\varphi}) 
           + b_E^h\bigl( (1 - \Pi^0) \vm{\varphi}, 
                   (1 - \Pi^0) \vm{\varphi} \bigr) \nonumber \\
        &= b_E^h(\widehat{\vm{M}}_0\vm{\Pi}_*^0,     
                 \widehat{\vm{M}}_0\vm{\Pi}_*^0 )   
            + \suku{  b_E^h\bigl( \vm{\varphi} ( \vm{I} - \vm{\Pi}^0) , 
                      \vm{\varphi} (\vm{I} - \vm{\Pi}^0) \bigr) } \nonumber \\       
        &= ( \vm{\Pi}_*^0)^T     
           b_E^h(\widehat{\vm{M}}_0 ,  \widehat{\vm{M}}_0 ) 
           \, \vm{\Pi}_*^0                  
           + (\vm{I} - \vm{\Pi}^0)^T \,  b_E^h(\vm{\varphi},
              \vm{\varphi}) \
             (\vm{I} - \vm{\Pi}^0 ) \nonumber \\
        &\approx ( \vm{\Pi}_*^0)^T     
           \vm{H}_0
           \, \vm{\Pi}_*^0                  
           + \rho |E| \, (\vm{I} - \vm{\Pi}^0)^T 
             (\vm{I} - \vm{\Pi}^0)  \nonumber \\
             &:= \vm{M}_E^c + \vm{M}_E^s,
\end{align}
where $\widehat{\vm{M}}_0 = \widehat{\vm{m}}_0 \otimes \vm{I}_d$,
$\vm{H}_0$ is a matrix of integrals of 
monomials up to order 2, which are listed in~\eqref{eq:integrals}, 
and we approximate $b_E^h(\vm{\varphi},\vm{\varphi})$ by
$\rho |E|$ in the stabilization term.

\subsubsection{Lumped mass matrix}
We perform mass lumping to construct diagonal mass matrices
using row-sum and diagonal scaling techniques.  In the
row-sum technique,\cite{Hughes:2000:FEM} the diagonal entries of the lumped mass matrix is given by 
\begin{equation}
M_{E,ii}^\ell = 
        \sum_{j=1}^{d N_E} M_{E,ij} \quad (d = 2,3).
\end{equation}
The diagonal entries of the lumped mass matrix using
the diagonal scaling procedure are given by~\cite{Hinton:1976:NML}
\begin{equation}
M_{E,ii}^\ell = 
        \dfrac{M_{E,ii}} {\texttt{trace} \, (\vm{M}_E) } d \rho |E| ,
\end{equation}
which ensures that mass conservation is satisfied ($\rho |E|$ is the mass in each direction).  Diagonal scaling always leads to
lumped mass matrices with positive diagonal entries, which is not guaranteed with the row-sum technique.
For convex polyhedral elements, we use the row-sum lumped mass technique, and for nonconvex polyhedral elements the diagonally-scaled lumped 
mass matrix is used.\cite{Park:2019:ONM,Park:2020:NRE}

\subsection{Estimation of critical time step} 
\label{subsec:timestep}
For stable explicit linear elastodynamic simulations, the
time step increment is subject to a restriction, which is given by the CFL condition:
\begin{equation}
  \Delta t \leq \frac{2}{\omega_\textrm{max}},
\end{equation}
where $\omega_{\textrm{max}}$ is the maximum natural frequency of vibration of the elastic solid. The natural frequency is related to the eigenvalue $\lambda$ of the generalized eigenproblem,
\begin{equation}
    \vm{K}\vm{d} = \lambda \vm{M} \vm{d},
\end{equation}
by the relation $\lambda = \omega^2 \ge 0$. By the element-eigenvalue inequality,~\cite{Fried:1972:BEE} we have 
$   \omega^* = \max \, \left (\omega_E^\textrm{max} \right) \ge \omega_\textrm{max} $ ,
and therefore a conservative lower bound
estimate for the stable critical time step is given by 
\begin{equation}
    \Delta t_\textrm{crit} = \frac{2}{\omega^*}
    \leq  \frac{2}{\omega_{max}} .
\end{equation}

\subsection{Homogeneous numerical integration scheme}\label{subsec:hni}
The homogeneous numerical integration (HNI) scheme~\cite{Chin:2015:NIH} allows to reduce integration of homogeneous functions
over arbitrary convex and nonconvex polytopes to integration over the
boundary facets of the polytope.
Let $f(\vx)$ be a positively homogeneous function of degree $q$ that is continuously differentiable:
\begin{equation}
f(\lambda \vx) = \lambda^q f(\vx) \quad (\lambda > 0),
\end{equation}
which satisfies Euler's homogeneous function theorem:
\begin{equation}\label{eq:euler}
q f(\vx) = \bigl( \nabla f(\vx), \vx \bigr) \quad \forall \vx \in \Re^d,
\end{equation}
where $(\cdot , \cdot )$ denotes the inner product of vectors in $\Re^d$
and $\vx$ is the position vector in $\Re^d$. 

Consider a polyhedral element $E$. The generalized Stokes's (divergence) theorem over $E$ can be stated as:
\begin{equation} \label{eq:div}
\int_E ( \nabla \cdot \vm{X} ) f ( \vm{x} ) \, d \vm{x} + 
\int_E \big( \vm{X}, \nabla f( \vm{x} ) \bigr) d \vm{x}
= \int_{\partial E} ( \vm{X}, \vm{n}) f( \vm{x} ) \, d \sigma ,
\end{equation}
where $d\sigma$ is the Lebesgue measure on $\partial E$.
For a homogeneous function $f$ and letting $\vm{X} := \vx$ be a vector
field and noting that $\nabla \cdot \vx = d$, 
\eqref{eq:div} reduces to
\begin{align} \label{eq:E}
\int_E f(\vx) \, d \vx = \frac{1}{d+q} 
\sum_{i=1}^m \frac{b_i}{\|\vm{a}_i\|} 
\suku{ \int_{\face_i} }  f(\vm{x}) \,d \sigma ,
\end{align}
where $\suku{ \partial E = \cup_{i=1}^m \face_i }$ 
consists of $m$ planar polygonal faces. Each \suku{$\face_i$} is a subset of the hyperplane ${\cal H}_i = \{ \vx : \vm{a}_i \cdot \vx = b_i \}$. 

On reapplying Stokes's theorem and Euler's homogeneous function theorem
on $\suku{\face_i \subset {\cal H}_i}$, we can write~\cite{Chin:2015:NIH}
\suku{
\begin{equation} \label{eq:Fi}
\int_{\face_i} f ( \vm{x} ) \, d \sigma =  \frac{1}{d + q - 1} \Biggl[  
\sum_{j \neq i}      
\int_{\face_{ij}} d_{ij} f ( \vm{x} ) \, d \nu 
+ \int_{\face_i} \big( \nabla f ( \vm{x} ), \vm{x}_0 \bigr) d \sigma \Biggr],
\end{equation}}
where $\vx_0 \in {\cal H}_i$ is an arbitrary point that 
satisfies $\vm{a}_i \cdot \vm{x}_0 = b_i$, 
\suku{$\face_{ij} \subset {\cal H}_{ij} = \face_i \cap \face_j$} is the $(d-2)$-dimensional facet, and
$d_{ij} := ( \vm{x} - \vm{x}_0 , \vm{n}_{ij} ) $ is the algebraic distance from $\vx_0$ to ${\cal H}_{ij}$. 

For the first-order VEM ($k = 1$), we require the computation of integrals over $E$ that involve monomials of maximum order 1 for the stiffness matrix and monomials of maximum order 2 for the mass matrix. We compute all such integrals by using a bottom-up approach, starting from the zeroth order monomial (volume computation) to order 2. In doing so, previously stored values are used to evaluate the second term (integral of monomials over the faces of the polyhedra) on the right-hand side of ~\eqref{eq:Fi}.  This leads to exact computation of these integrals.\cite{Chin:2020:AEM} Referring 
to~\eqref{eq:xietazeta}, we use HNI to compute the following
integrals: 
\begin{align} \label{eq:integrals}
\int_E d \vx = |E|, \ \ 
\int_E \xi \, d \vx &= 
\frac{1}{h_E} \left[ \int_E x \, d \vx  - x_E |E| \right] , \ \
\int_E \eta \, d \vx = 
\frac{1}{h_E} \left[ \int_E y \, d \vx  - y_E |E| \right] , \ \
\int_E \zeta \, d \vx = 
\frac{1}{h_E} \left[ \int_E z \, d \vx  - z_E |E| \right] , \nonumber
\\
\int_E \xi \eta \, d \vx &= 
\frac{1}{h_E^2} \left[ \int_E xy \, d \vx  - y_E \int_E x \, d\vx 
-x_E \int_E y d\vx + x_E y_E |E| \right] , \nonumber \\
\int_E \xi \zeta \, d \vx &= 
\frac{1}{h_E^2} \left[ \int_E xz \, d \vx  - z_E \int_E x \, d\vx 
-x_E \int_E z d\vx + x_E z_E |E| \right] , \nonumber \\
\int_E \eta \zeta \, d \vx &= 
\frac{1}{h_E^2} \left[ \int_E yz \, d \vx  - z_E \int_E y \, d\vx 
- y_E \int_E z d\vx + y_E z_E |E| \right] , \\
\int_E \xi^2 \, d \vx &= 
\frac{1}{h_E^2} \left[ \int_E x^2 \, d \vx  - 2x_E \int_E x \, d\vx 
                       + x_E^2 |E| \right] , \nonumber \\
\int_E \eta^2 \, d \vx &= 
\frac{1}{h_E^2} \left[ \int_E y^2 \, d \vx  - 2y_E \int_E y \, d\vx 
                       + y_E^2 |E| \right] , \nonumber \\     
\int_E \zeta^2 \, d \vx &= 
\frac{1}{h_E^2} \left[ \int_E z^2 \, d \vx  - 2z_E \int_E z \, d\vx 
                       + z_E^2 |E| \right] \nonumber ,                          
\end{align}
where $|E|$ is the volume of the polyhedral element $E$.

\section{Numerical Results}\label{sec:results}
We use the element-eigenvalue inequality to obtain an estimate of the
critical time step for linear elastodynamic simulations. FEM is executed on Delaunay meshes and VEM is adopted on polyhedra that are formed via agglomeration of tetrahedral and prismatic elements. In all computations, material properties of steel are used:
$E_Y = 210$ GPa, $\nu = 0.3$ and $\rho = 7800$ kg/$\textrm{m}^3$.

\subsection{Two-dimensional study}
\begin{figure}[!h]
  \centering
  \begin{subfigure}{0.48\textwidth}
      \includegraphics[width=0.78\textwidth]{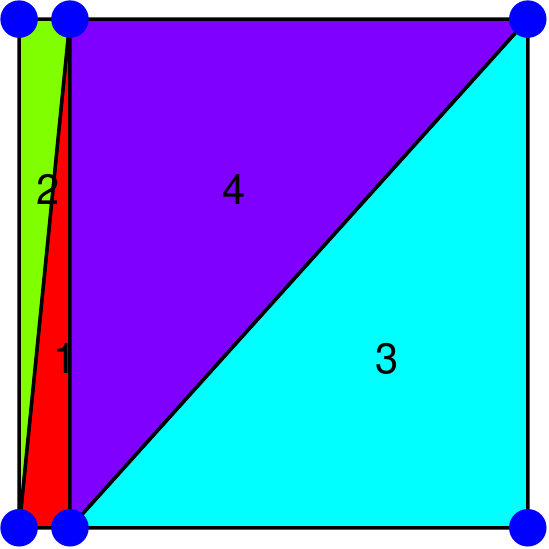}
      \caption{FEM ($\epsilon = 10^{-1}$)}\label{fig:meshes2D-a}
  \end{subfigure} \hfill
  \begin{subfigure}{0.48\textwidth}
      \includegraphics[width=0.78\textwidth]{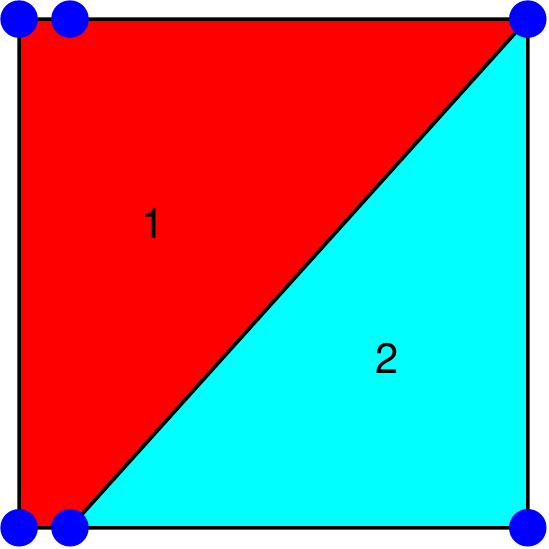}
      \caption{VEM ($\epsilon = 10^{-1}$)}\label{fig:meshes2D-b}
  \end{subfigure}
  \begin{subfigure}{0.48\textwidth}
      \includegraphics[width=0.78\textwidth]{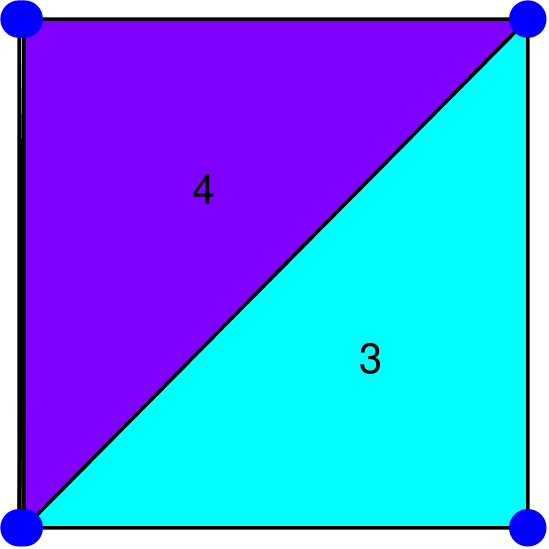}
      \caption{FEM ($\epsilon = 10^{-2}$)}\label{fig:meshes2D-c}
  \end{subfigure} \hfill
  \begin{subfigure}{0.48\textwidth}
      \includegraphics[width=0.78\textwidth]{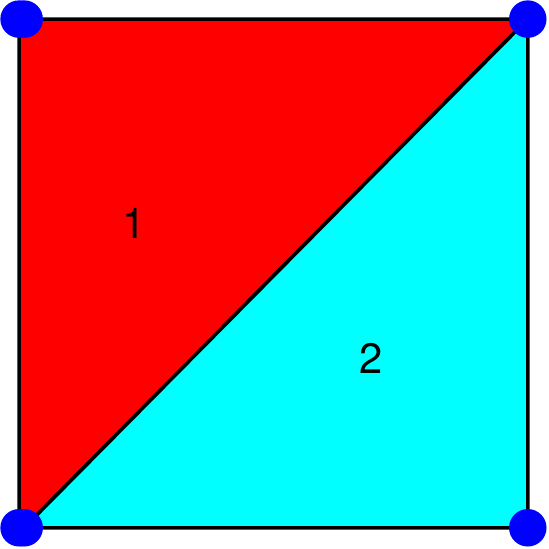}
      \caption{VEM ($\epsilon = 10^{-2}$)}\label{fig:meshes2D-d}
  \end{subfigure}
  \begin{subfigure}{0.48\textwidth}
      \includegraphics[width=0.8\textwidth]{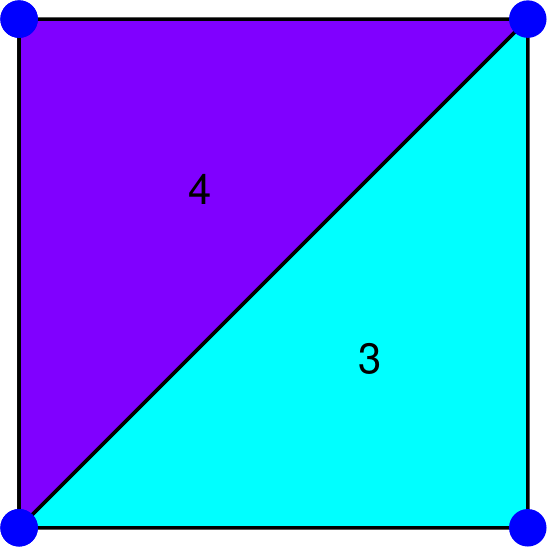}
      \caption{FEM ($\epsilon = 10^{-5}$)}\label{fig:meshes2D-e}
  \end{subfigure}\hfill
  \begin{subfigure}{0.48\textwidth}
      \includegraphics[width=0.8\textwidth]{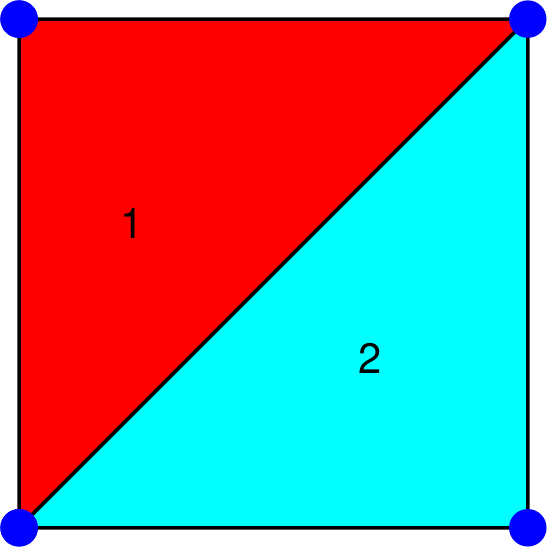}
      \caption{VEM ($\epsilon = 10^{-5}$)}\label{fig:meshes2D-f}
  \end{subfigure}
  \caption{Triangular finite element and polygonal meshes that are
           used 
           to determine the critical time step. 
           }\label{fig:meshes2D}
\end{figure}

We assume plane strain conditions and
consider sequence of finite
element meshes in which the length of the smallest edge is monotonically
decreased ($\epsilon = 10^{-1}, \, 10^{-2}, \,10^{-5}, \, 10^{-8}$). 
Three meshes from this sequence for FE and VEM (polygonal meshes) are
presented in~\fref{fig:meshes2D}. 
The results for the maximum natural frequency as well
as the improvements in $\Delta t_{\textrm{crit}}$ that are
obtained using VEM are listed
in~\tref{tab:eig}. The second element serves as a reference since this element is
well-shaped (right-angled nearly isosceles triangle). The overall maximum frequency
on the FE meshes increases substantially
from $ 1.5 \times 10^5$ Hz to $ 1.5 \times 10^{12}$ Hz, whereas for the VEM the increase is from
$ 3.6 \times 10^4$ Hz to $ 5.1 \times 10^4$ Hz. For the VEM, the maximum frequency
is relatively insensitive to the minimum edge length.
On using the VEM on agglomerated element much larger time steps can be used in comparison to the FEM, with a factor of $4$
when $\epsilon = 1$ to a factor of $10^7$ when
$\epsilon = 10^{-8}$.
\begin{table}
\centering
\caption{Critical time step in 2D for FEM
         (Delaunay meshes) and 
         VEM using lumped (row-sum) mass matrix. For the VEM, bad triangles (interior angle close to $0^\circ$) are
          agglomerated with neighboring triangles to form a polygonal element.}
          \label{tab:eig}
\begin{tabular}{ |c | c | c   c | c   c | c | }
\hline
&&&&&& \\
Mesh & Method & $E$ & $\omega_{\textrm{max}}^E$ & $E$ & $\omega_{\textrm{max}}^E$ & $\dfrac{\Delta t_{\textrm{crit}}^{\textrm{VEM}}} {\Delta t_{\textrm{crit}}^{\textrm{FEM}}}$ \\
&&&&&& \\ \hline
& & & & & & \\ 
\ref{fig:meshes2D-a} ($\epsilon = 10^{-1}$) & FEM & 1 & $1.5 \times 10^{5}$ & 3 & $2.0 \times 10^{4}$ & $\sim 4$ \\
\ref{fig:meshes2D-b} ($\epsilon = 10^{-1}$) & VEM & 1 & $3.6 \times 10^{4}$ & 2 & $2.0 \times 10^{4}$ & \\
& & & & & & \\
\ref{fig:meshes2D-c} ($\epsilon = 10^{-2}$) & FEM & 1 & $1.5 \times 10^{6}$ & 3 & $1.9 \times 10^{4}$ & $\sim \suku{30} $ \\
\ref{fig:meshes2D-d} ($\epsilon = 10^{-2}$) & VEM & 1 & $4.9 \times 10^{4}$ & 2 & $1.9 \times 10^{4}$ & \\
& & & & & & \\
\ref{fig:meshes2D-e} ($\epsilon = 10^{-5}$) & FEM & 1 & $1.5 \times 10^{9}$ & 3 & $1.9 \times 10^{4}$ & $\sim 10^4$ \\
\ref{fig:meshes2D-f} ($\epsilon = 10^{-5}$) & VEM & 1 & $5.1 \times 10^{4}$ & 2 & $1.9 \times 10^{4}$ & \\
& & & & & & \\
$\epsilon = 10^{-8}$ & FEM & 1 & $1.5 \times 10^{12}$ & 3 & $1.9 \times 10^{4}$ & $\sim 10^7$ \\
$\epsilon = 10^{-8}$ & VEM & 1 & $5.1 \times 10^{4}$ & 2 & $1.9 \times 10^{4}$ & \\
& & & & & &  \\ \hline
\end{tabular}
\end{table}

\subsection{Three-dimensional study}
We consider one case with prismatic elements and then three cases that contain poor-quality tetrahedral elements. For the prismatic element, we use the row-sum technique to form the lumped mass matrix, whereas for the polyhedral elements (in general, nonconvex polyhedra) that are formed by agglomerating sliver (kite) and spire tetrahedra,
diagonal scaling is used to form the lumped mass matrix. 
\suku{In the 3D study, since we use meshes with elements whose diameters are ${\cal O}(1)$, we 
set $\alpha_0 = 1$
in the stabilization term of the stiffness matrix.}

\subsubsection{Prismatic elements} 
\begin{figure}
  \centering
  \begin{subfigure}{0.49\textwidth}
      \includegraphics[width=0.85\textwidth]{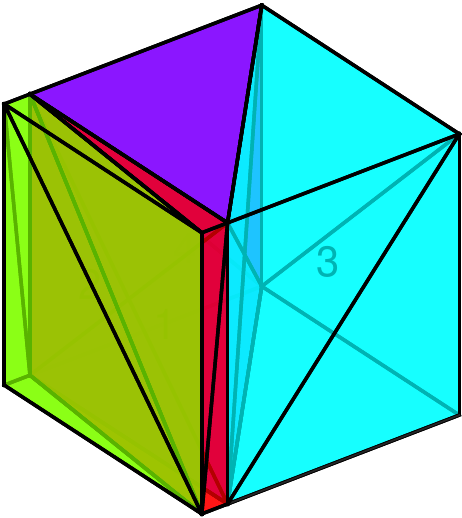}
      \caption{FEM ($\epsilon = 10^{-1}$)}\label{fig:meshes3D-a}
  \end{subfigure} \hspace*{0.06in}
  \begin{subfigure}{0.49\textwidth}
      \includegraphics[width=0.85\textwidth]{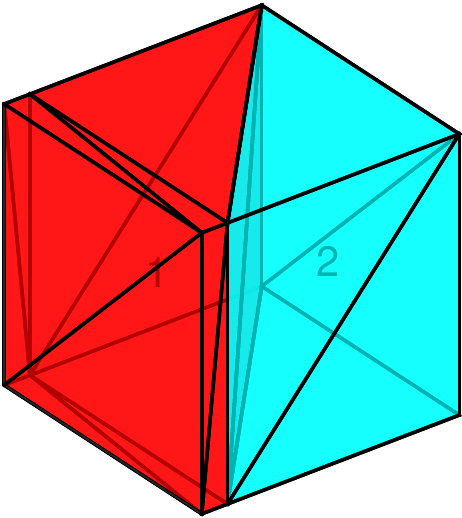}
      \caption{VEM ($\epsilon = 10^{-1}$)}\label{fig:meshes3D-b}
  \end{subfigure}
  \begin{subfigure}{0.49\textwidth}
      \includegraphics[width=0.85\textwidth]{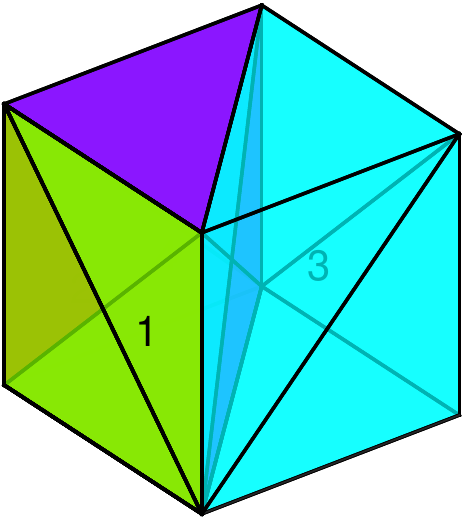}
      \caption{FEM ($\epsilon = 10^{-5}$)}\label{fig:meshes3D-c}
  \end{subfigure} \hspace*{0.06in}
  \begin{subfigure}{0.49\textwidth}
      \includegraphics[width=0.85\textwidth]{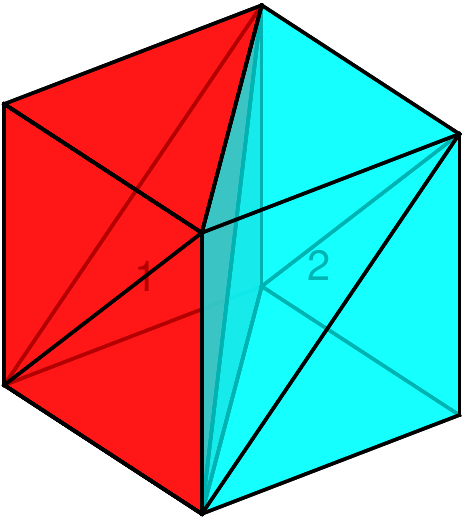}
      \caption{VEM ($\epsilon = 10^{-5}$)}\label{fig:meshes3D-d}
  \end{subfigure}
  \caption{Prismatic meshes (FEM) and polyhedral meshes (triangular faces) used in the study
           to determine the critical time step.
           }\label{fig:meshes3D}
\end{figure}

The mesh that was considered in two dimensions is extended to three dimensions by extruding the two-dimensional meshes in the
third dimension, which results in prismatic elements.  The agglomeration is done
in a similar manner as in 2D. Two polyhedral virtual elements are created and the
polygonal faces are represented as the union of triangular faces.
Meshes for FE and VEM are presented in~\fref{fig:meshes3D}, and  the
results for the maximum natural frequency as well
as the ratio of $\Delta t_{\textrm{crit}}$ obtained using VEM and FEM are listed in~\tref{tab:eig3D}.
The trends are similar to the 2D case. The faces with very small areas do not
affect the maximum frequency in the virtual element computations.
\begin{table}
\centering
\caption{Critical time step using lumped (row-sum)
         mass matrix for stretched prismatic finite elements and agglomerated 
         virtual elements.
         For the virtual element computations, prismatic elements with vanishing
         face areas are agglomerated with their neighbors to create a hexahedral element
         with many triangular faces.}
         \label{tab:eig3D}
\begin{tabular}{ |c | c | c   c | c   c | c |}
\hline
&&&&& & \\
Mesh & Method & $E$ & $\omega_{\textrm{max}}^E$ (Hz) & $E$ & $\omega_{\textrm{max}}^E$ (Hz) & 
$\dfrac{\Delta t_{\textrm{crit}}^{\textrm{VEM}}} {\Delta t_{\textrm{crit}}^{\textrm{FEM}}}$ \\
&&&&&& \\ \hline
& & & & & & \\ 
\ref{fig:meshes3D-a} ($\epsilon = 10^{-1}$) & FEM & 1 & $1.7 \times 10^{5}$ & 3 & $2.6 \times 10^{4}$ & $\sim 3$ \\
\ref{fig:meshes3D-b} ($\epsilon = 10^{-1}$) & VEM & 1 & $4.9 \times 10^{4}$ & 2 & $2.6 \times 10^{4}$ &  \\
&&&&&& \\
\ref{fig:meshes3D-c} ($\epsilon = 10^{-5}$) & FEM & 1 & $1.7 \times  10^{9}$ & 3 & $2.5 \times 10^{4}$ & $\sim 10^4$ \\
\ref{fig:meshes3D-d} ($\epsilon = 10^{-5}$) & VEM & 1 & $1.0  \times 10^{5}$ & 2 & $2.5 \times 10^{4}$ & \\
& & & & & & \\ \hline
\end{tabular}
 \end{table}

\subsubsection{Tetrahedral elements} 
Tetrahedral meshes with wedge, kite configuration and spires are considered.\cite{Cheng:2000:SE} Agglomeration of a badly-shaped tetrahedral element with a neighboring element
is performed to construct a polyhedral virtual element that is used
in the computations.

\begin{figure}
  \centering
  \begin{subfigure}{0.45\textwidth}
      \includegraphics[width=\textwidth]{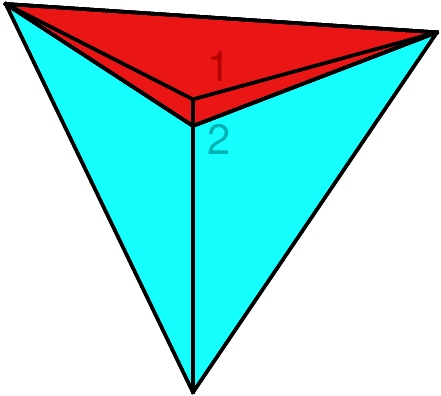}
      \caption{FEM ($\epsilon = 10^{-1}$)}\label{fig:tet-wedge-a}
  \end{subfigure}
  \begin{subfigure}{0.45\textwidth}
      \includegraphics[width=\textwidth]{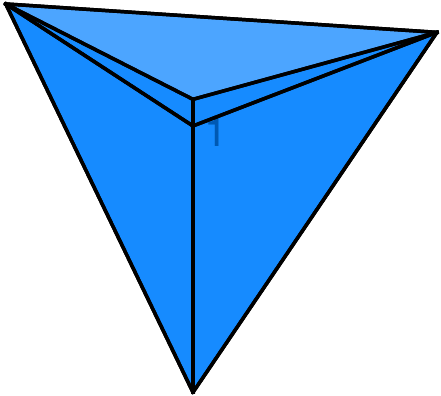}
      \caption{VEM ($\epsilon = 10^{-1}$)}\label{fig:tet-wedge-b}
  \end{subfigure}
  \caption{Influence of a wedge tetrahedral element on the critical time step. 
  (a) Two tetrahedral finite elements with one of them being a wedge. 
  (b) Agglomeration of the two tetrahedral elements into a single 
      polyhedral virtual element.}\label{fig:tet-wedge}
\end{figure}

\smallskip
\noindent {\bf Wedge tetrahedron.}
Consider the two-element tetrahedral configuration shown
in~\fref{fig:tet-wedge-a}. The tetrahedral element above the $xy$-plane is a wedge tetrahedron (one dihedral angle is close to $0^\circ$) with $\epsilon$ being the $z$-coordinate of the node above the $xy$-plane. 
The tetrahedral element below the $xy$-plane is well-shaped. 
The meshes shown in~\fref{fig:tet-wedge} are for 
$\epsilon = 10^{-1}$. 
The two elements are combined to form a single polyhedral element with six faces (see~\fref{fig:tet-wedge-b}).
Numerical results for the maximum natural frequency for
$\epsilon = 10^{-1},\,10^{-3},\,10^{-5}$ and the ratio of
the critical time step for VEM to that obtained using the
FEM are listed
in~\tref{tab:eig3D_tet_wedge}. For $\epsilon = 10^{-5}$,
the ratio of the critical time step of VEM to that of
FEM is ${\cal O}(10^4)$.
\begin{table}
\centering
\caption{Critical time step using 
         lumped (diagonal scaling) mass matrix for 
         tetrahedral finite elements with a wedge and 
         an agglomerated virtual element. The
         agglomerated virtual element has six triangular
         faces.}
         \label{tab:eig3D_tet_wedge}
\begin{tabular}{ |c | c | c   c | c   c | c |}
\hline
&&&&& & \\
$\epsilon$ & Method & $E$ & $\omega_{\textrm{max}}^E$ (Hz) & $E$ & $\omega_{\textrm{max}}^E$ (Hz) & 
$\dfrac{\Delta t_{\textrm{crit}}^{\textrm{VEM}}} {\Delta t_{\textrm{crit}}^{\textrm{FEM}}}$ \\
&&&&&& \\ \hline
& & & & & & \\ 
& FEM & 1 & $1.7 \times 10^{5}$ & 2 & $2.5 \times 10^{4}$ &  \\
$10^{-1}$ & & & & & & $\sim 3$ \\
& VEM & 1 & $4.3 \times 10^{4}$ & -- & -- &  \\
&&&&&& \\ &&&&&& \\
& & & & & & \\ 
& FEM & 1 & $1.7 \times 10^{7}$ & 2 & $2.5 \times 10^{4}$ &  \\
$10^{-3}$ & & & & & & $\sim 10^2$ \\
& VEM & 1 & $4.6 \times 10^{4}$ & -- & -- &  \\
&&&&&& \\ &&&&&& \\
& FEM & 1 & $1.7 \times 10^{9}$ & 2 & $2.5 \times 10^{4}$ &  \\
$10^{-5}$ & & & & & & $\sim 10^4$ \\
& VEM & 1 & $4.6 \times 10^{4}$ & -- & -- &  \\
&&&&&& \\ \hline
\end{tabular}
 \end{table}

\begin{figure}
  \centering
  \begin{subfigure}{0.45\textwidth}
      \includegraphics[width=\textwidth]{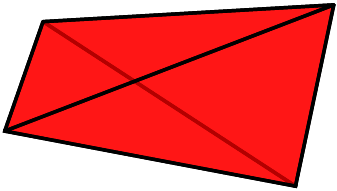}
      \caption{FEM ($\epsilon = 10^{-1}$)}\label{fig:tet-kite-a}
  \end{subfigure}
  \begin{subfigure}{0.45\textwidth}
      \includegraphics[width=\textwidth]{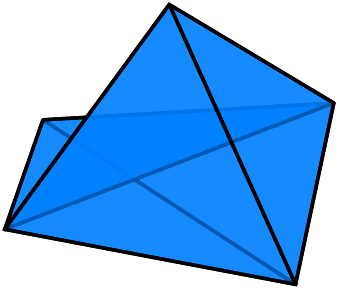}
      \caption{VEM ($\epsilon = 10^{-1}$)}\label{fig:tet-kite-b}
  \end{subfigure}
  \caption{Influence of tetrahedral kite element on 
           the critical time step. 
           (a) Tetrahedral kite element and
           (b) Agglomeration of the kite element with 
            a neighboring tetrahedral element to form 
            a polyhedral virtual element.}
           \label{fig:tet-kite}
\end{figure}

\smallskip
\noindent {\bf Sliver (kite) tetrahedron.}
We consider a sliver (kite) tetrahedron element. The coordinates of the nodes of this tetrahedron are:
$(-1,0,\epsilon)$, $(1,0,\epsilon)$, $(0,-1,-\epsilon)$, 
and $(0,1,\epsilon)$. Two dihedral angles of this tetrahedron are close to $180^\circ$, and the other four dihedral angles are close
to $0^\circ$. We join this tetrahedron to a neighboring 
tetrahedron with the additional vertex $(0,0,1$) to form
a polyhedral (nonconvex) virtual element.  The kite tetrahedron and the
polyhedral virtual element are shown in~\fref{fig:tet-kite}.

Numerical results for the maximum natural frequency for
$\epsilon = 10^{-1},\,10^{-5}$, and the ratio of
the critical time step for VEM to that obtained using the 
FEM are listed
in~\tref{tab:eig3D_tet_kite}. We observe that for
$\epsilon = 10^{-5}$, the critical
time step estimate for the VEM 
polyhedral element is four orders greater than that of FEM.
As $\epsilon$ decrease this ratio increases; the critical step in the VEM is insensitive to
$\epsilon$. As $\epsilon \rightarrow 0$, the stable time
step for the VEM is of the same order as a well-shaped tetrahedral finite element.
\begin{table}
\centering
\caption{Critical time step using lumped
         (diagonal scaling) mass matrix for a tetrahedral  
         kite finite element and an agglomeration of it with a
         neighboring elements into a polyhedral (nonconvex)
         virtual element. The
         agglomerated virtual element has six triangular
         faces.}
         \label{tab:eig3D_tet_kite}
\begin{tabular}{ |c | c  c | c  c | c | }
\hline
&&&&&  \\
$\epsilon$ & Method & $\omega_{\textrm{max}}^{\textrm{FEM}}$ (Hz) & 
Method & $\omega_{\textrm{max}}^{\textrm{VEM}}$ (Hz) & 
$\dfrac{\Delta t_{\textrm{crit}}^{\textrm{VEM}}} {\Delta t_{\textrm{crit}}^{\textrm{FEM}}}$ \\
&&&&& \\ \hline
& & & & & \\ 
$10^{-1}$ & FEM & $6.0 \times 10^{4}$ & VEM & $3.1 \times 10^{4}$ &  $\sim 2$ \\
& & & & & \\ 
$10^{-5}$ & FEM & $6.0 \times 10^{8}$ & VEM & $5.2 \times 10^{4}$ &  $\sim 10^4$ \\
&&&&& \\ \hline
\end{tabular}
\end{table}

\begin{figure}[!h]
  \centering
  \begin{subfigure}{0.58\textwidth}
      \includegraphics[width=0.9\textwidth]{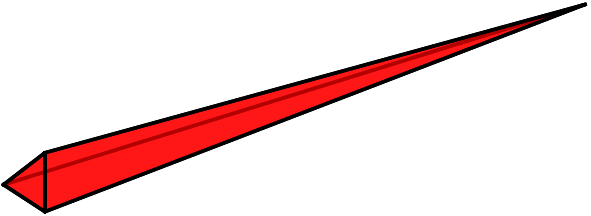}
      \caption{FEM ($\epsilon = 10^{-1}$)}\label{fig:tet-spire-a}
  \end{subfigure}
  \begin{subfigure}{0.38\textwidth}
      \includegraphics[width=0.9\textwidth]{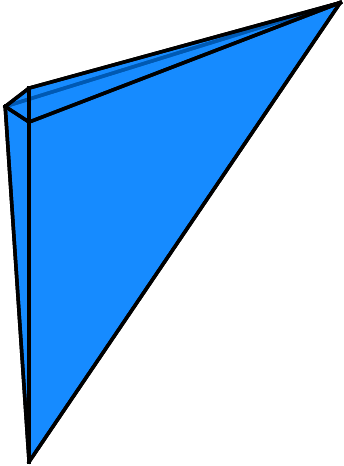}
      \caption{VEM, Case $A$ ($\epsilon = 10^{-1}$)}\label{fig:tet-spire-b}
  \end{subfigure}
  
  \vspace*{0.6in}
  \begin{subfigure}{0.48\textwidth}
     \hspace*{0.2in}
      \includegraphics[width=0.8\textwidth]{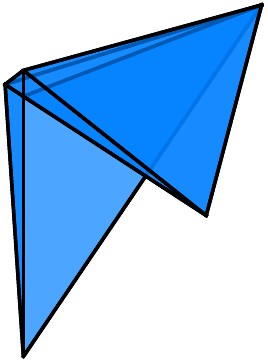}
      \caption{VEM, Case $B$ ($\epsilon = 10^{-1}$)}\label{fig:tet-spire-c}
  \end{subfigure}\hspace*{0.4in}
  \begin{subfigure}{0.48\textwidth}
      \includegraphics[width=0.8\textwidth]{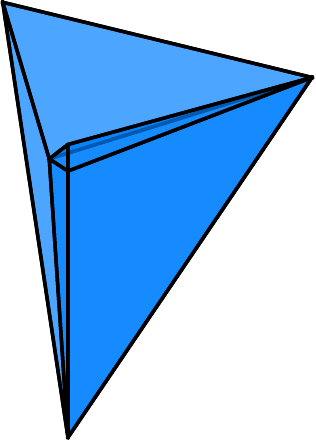}
      \caption{VEM, Case $C$ ($\epsilon = 10^{-1}$)}\label{fig:tet-spire-d}
  \end{subfigure}
  \caption{Influence of tetrahedral spire element on the 
            critical time step. 
           (a) Tetrahedral spire element. Agglomeration of the
           spire with 
           (b) one neighboring element
           (case $A$), and (c), (d) two neighboring elements
           (cases $B$ and $C$) to form a polyhedral virtual
           element.}
           \label{fig:tet-spire}
\end{figure}

\smallskip
\noindent {\bf Spire tetrahedron.}
We consider a spire tetrahedron, which is a tetrahedron with one tiny face and three long edges. The nodal coordinates
of this tetrahedron are chosen as: 
$(0,0,0)$, $(0,\epsilon,0)$, $(0,0,\epsilon)$, 
and $(1,0,0)$. For the neighboring tetrahedron, we consider the
vertices $(0,0,-1)$, $(0,-1,0)$ and $(1/2,1,0)$. We join the spire tetrahedron with one or two neighboring tetrahedra to form a polyhedral
virtual element. Three such cases are considered that we label as
$A$, $B$, and $C$. For $\epsilon = 10^{-1}$, the spire tetrahedron and the polyhedral virtual element (cases $A$,
$B$, and $C$) that
are formed via agglomeration are shown in~\fref{fig:tet-spire}.

Numerical results for the maximum natural frequency for
$\epsilon = 10^{-1},\,10^{-5}$, and the ratio of
the critical time step for VEM (cases $A$, $B$ and $C$) vis-{\`a}-vis FEM are listed
in~\tref{tab:eig3D_tet_spire}. For cases $A$ and $B$ with the
virtual element (severly distored nonconvex elements), the critical time step is appreciably greater than that for the FEM. For case $C$ (nearly convex element), the VEM outperforms the FEM, and the trends are similar to results for other problems that are reported earlier. For
$\epsilon = 10^{-5}$, the critical
time step estimate for the VEM 
polyhedral element (case $C$) is four orders greater than that of FEM. We find that as $\epsilon \rightarrow 0$, the stable time
step for this case is of the same order as a well-shaped tetrahedral finite element. For instance, if $\epsilon = 10^{-8}$, we find that $\omega_{\textrm{max}}^{\textrm{VEM}} = 5.9 \times 10^5$ Hz. To promote
better conditioning of the stiffness matrix, it is preferable to
agglomerate with more neighbors to produce a polyhedral virtual element with nonvanishing volume.
\begin{table}
\centering
\caption{Critical time step using lumped
         (diagonal scaling) mass matrix for a tetrahedral  
         spire finte element and an agglomeration of it with
         neighboring elements into a polyhedral (nonconvex)
         virtual element.}
         \label{tab:eig3D_tet_spire}
\begin{tabular}{ |c | c  c | c  c | c | }
\hline
&&&&&  \\
$\epsilon$ & Method & $\omega_{\textrm{max}}^{\textrm{FEM}}$ (Hz) & 
Method & $\omega_{\textrm{max}}^{\textrm{VEM}}$ (Hz) & 
$\dfrac{\Delta t_{\textrm{crit}}^{\textrm{VEM}}} {\Delta t_{\textrm{crit}}^{\textrm{FEM}}}$ \\
&&&&& \\ \hline
& & & & & \\ 
& & & VEM (Case A)& $3.2 \times 10^{5}$ &  $\sim 1$ \\
$10^{-1}$ & FEM & $2.2 \times 10^{5}$ & VEM (Case B)& $3.9 \times 10^{5}$ &  $\sim 1$ \\
& & & VEM (Case C)& $5.5 \times 10^{4}$ &  $\sim 6$ \\
& & & & & \\ 
& & & VEM (Case A)& $3.3 \times 10^{9}$ &  $\sim 1$ \\
$10^{-5}$ & FEM & $2.2 \times 10^9$ & VEM (Case B)& $2.2 \times 10^{13}$ &  $ \sim 10^{-4}$ \\
& & & VEM (Case C)& $5.9 \times 10^{4}$ &  $\sim 10^4$ \\
&&&&& \\ \hline
\end{tabular}
\end{table}

\suku{
\subsection{Explicit dynamics simulations on a 
tapered beam}\label{subsec:explicit_dynamics}
To establish the robustness and efficiency of meshes with agglomerated
polyhedral virtual elements over poor-quality 
tetrahedral finite element meshes, we present elastodynamic simulations
for a manufactured problem. To this end, we begin with a two-dimensional beam of dimensions $[0,4] \times [0,1]$ that is meshed with bilinear finite elements.  We cut this mesh with an inclined line that results in some elements having
poor quality. For FEM, all 2D elements are split into triangles. For VEM, there are uncut hexahedral elements as well as polygonal (greater than four edges) elements that are formed as agglomerations of the triangles. Sample two-dimensional 
meshes of the tapered beam are shown in~\fref{fig:mesh_AB}. These meshes are extruded (one element thickness) in the $z$-direction. The 2D triangular elements yield prismatic elements in 3D. Each prism is meshed into three tetrahedral elements. This results in poor-quality tetrahedra in some of the regions in the vicinity of the cut plane. For the virtual element analysis, the elements are hexahedra away from the cutting plane.  Near the cutting plane the 2D triangles are aggregated into polygons as shown, then extruded in $z$ to form polyhedral virtual elements. Every rectangular face of the 3D extrusions is partitioned into two triangles. Each hexahedral virtual element has 12 triangular faces, whereas each agglomerated virtual element has 20 triangular faces. 

Three-dimensional meshes are generated by extruding
the two-dimensional meshes shown in~\fref{fig:mesh_AB}. We refer to case $A$
as the meshes shown in~\fref{fig:mesh_AB-a} and~\fref{fig:mesh_AB-b},
whereas case $B$ refers to the meshes shown in~\fref{fig:mesh_AB-c} and~\fref{fig:mesh_AB-d}.  Though not visible at the scale shown, the cut line is distinct in each case, with
case $B$ designed to result in sliver tetrahedra of poorer quality when compared to
case $A$. For cases $A$ and $B$, the number of nodes and elements for the tetrahedral meshes (and
likewise for the polyhedral meshes) are the same. The
number of nodes in all meshes is 1282.
The number of elements in the tetrahedral finite element 
mesh is 3456 and that in the polyhedral virtual element mesh is 549.  
The estimates for the maximum natural frequency in cases $A$ and $B$ are listed
in~\tref{tab:explicit_dynamics}. In~\fref{fig:ux_vs_time-b} (case $B$), 
since $\omega_{\textrm{max}}^{FEM} = {\cal O}(10^{16})$, 
the global system (stiffness and lumped mass) matrices
($\vm{K}\vm{d} = \omega^2 \vm{M}\vm{d})$ are used to compute
the maximum frequency in the FE simulations (see~\tref{tab:explicit_dynamics}).
For cases $A$ and $B$, the critical time step in the VEM is 275 and $4.5 \times 10^4$ times larger than the critical time step in the FEM, respectively.
\begin{figure}
  \centering
  \begin{subfigure}{0.46\textwidth}
      \includegraphics[width=\textwidth]{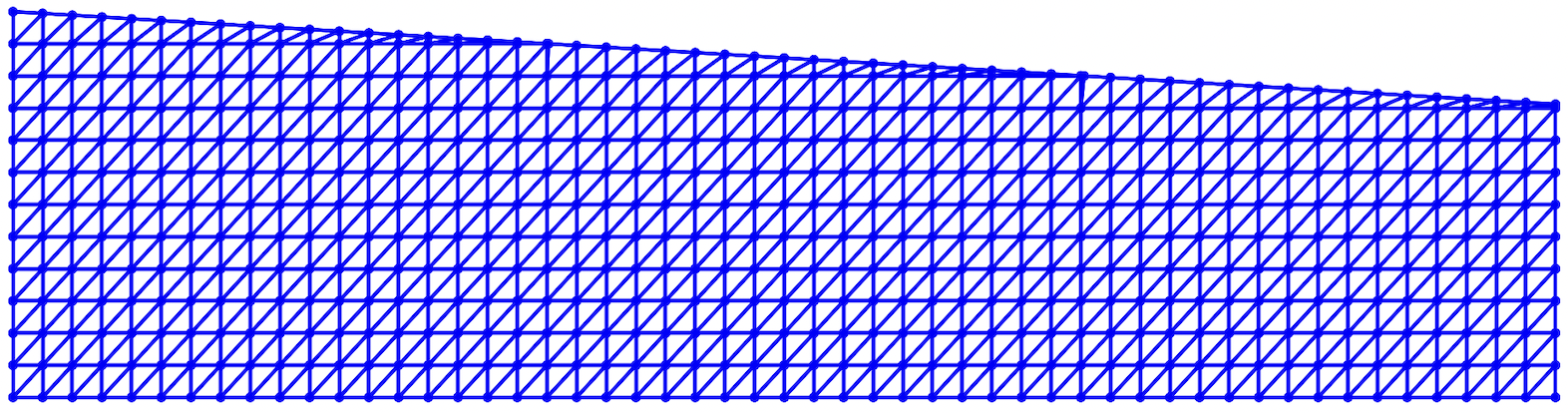}
      \caption{} \label{fig:mesh_AB-a}
  \end{subfigure} \hfill
  \begin{subfigure}{0.46\textwidth}
      \includegraphics[width=\textwidth]{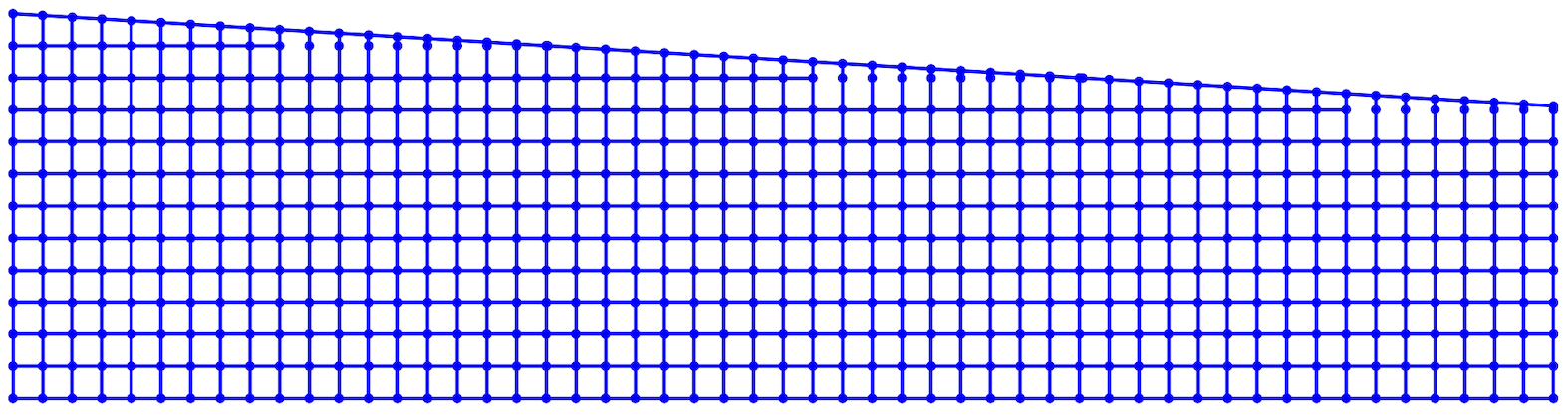}
      \caption{} \label{fig:mesh_AB-b}
  \end{subfigure} \hfill
  \begin{subfigure}{0.46\textwidth}
     \includegraphics[width=\textwidth]{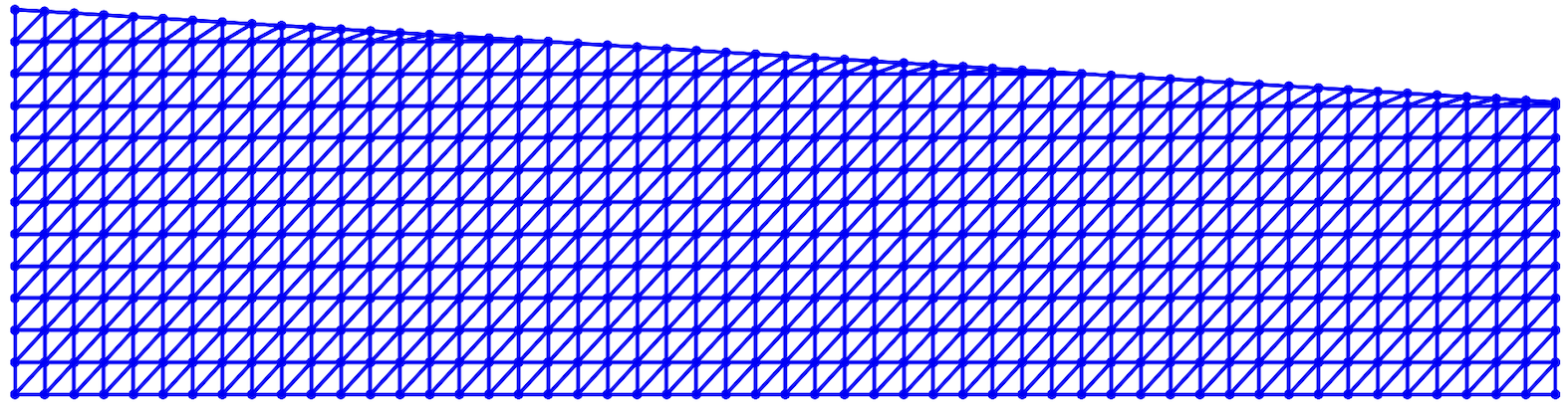}
     \caption{} \label{fig:mesh_AB-c}
  \end{subfigure} \hfill
  \begin{subfigure}{0.46\textwidth}
     \includegraphics[width=\textwidth]{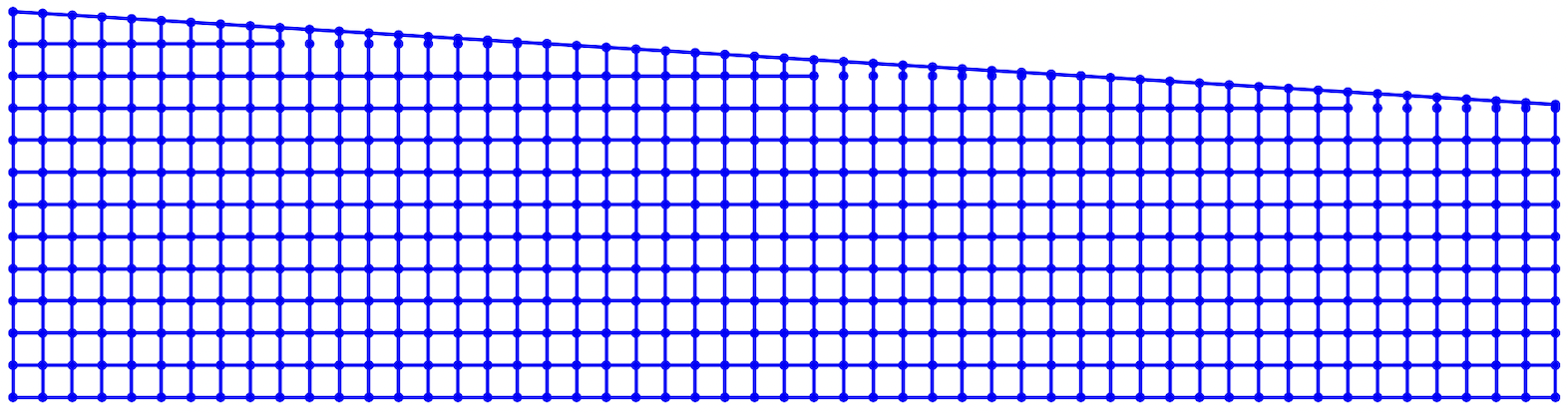}
     \caption{} \label{fig:mesh_AB-d}
  \end{subfigure}
  \caption{Two-dimensional meshes of the tapered beam, which are 
  constructed by cutting a structured mesh by a line.
  (a), (b) Triangular and polygonal meshes (Case $A$); and
  (c), (d) Triangular and polygonal meshes (Case $B$).}
           \label{fig:mesh_AB}
\end{figure}
\begin{table}[!hbt]
\centering
\caption{Critical time step for FEM and VEM
         in the explicit dynamics simulations.}
         \label{tab:explicit_dynamics}
\begin{tabular}{ | c |  c  c  c | }
\hline
&&&  \\
Meshes & $\omega_{\textrm{max}}^{\textrm{FEM}}$ (Hz)  
& $\omega_{\textrm{max}}^{\textrm{VEM}}$ (Hz) & 
$\dfrac{\Delta t_{\textrm{crit}}^{\textrm{VEM}}} {\Delta t_{\textrm{crit}}^{\textrm{FEM}}}$
\\
&&& \\ \hline
& & & \\ 
Case A & $8.0 \times 10^{7}$ & 
$2.9 \times 10^{5}$ &  $275$ \\
Case B & $1.3 \times 10^{10}$\textsuperscript{*}  & 
$2.9 \times 10^{5}$ &  $4.5 \times 10^4$ \\
&&& \\ \hline
\multicolumn{4}{l}{\textsuperscript{*}\footnotesize{Maximum natural
frequency is based on the global elastodynamic eigenproblem}}
\end{tabular}
\end{table}

In~\fref{fig:3Dmesh_A}, the surface triangular elements for case $A$
are depicted along with the 3D polyhedral virtual elements that have
$20$ faces. 
The tetrahedron with the worst condition number and the agglomerated polyhedron virtual element (20 faces) that contains the vertices of this tetrahedron are 
shown in~Figures~\ref{fig:tetpoly_mesh_A-a} and~\ref{fig:tetpoly_mesh_A-b},
respectively. 
\begin{figure}
  \centering
  \begin{subfigure}{0.46\textwidth}
     \includegraphics[width=\textwidth]{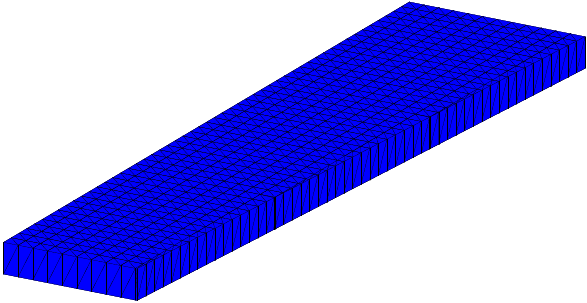}
     \caption{}  \label{fig:3Dmesh_A-a}
  \end{subfigure} \hfill
  \begin{subfigure}{0.46\textwidth}
      \includegraphics[width=\textwidth]{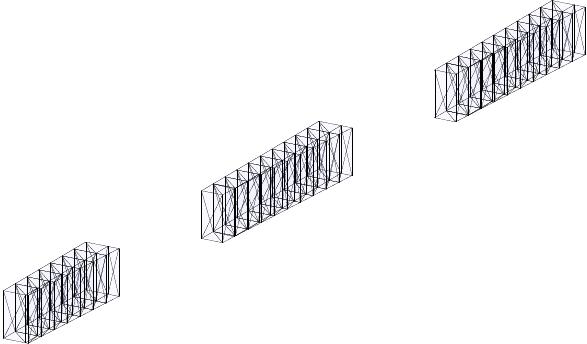}
     \caption{}  \label{fig:3Dmesh_A-b}
   \end{subfigure}
  \caption{(a) Surface triangular elements for
  tetrahedral and polyhedral virtual element meshes and
  (b) Agglomerated polyhedral virtual elements that have $20$ faces.}
           \label{fig:3Dmesh_A}
\end{figure}
\begin{figure}
  \centering
  \begin{subfigure}{0.44\textwidth}
     \includegraphics[height=2in,width=\textwidth]{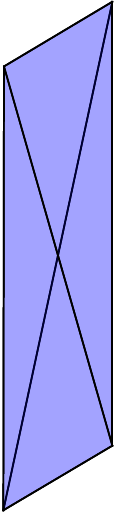}
     \caption{}  \label{fig:tetpoly_mesh_A-a}
  \end{subfigure} \hfill
  \begin{subfigure}{0.44\textwidth}
      \includegraphics[height=2in,width=\textwidth]{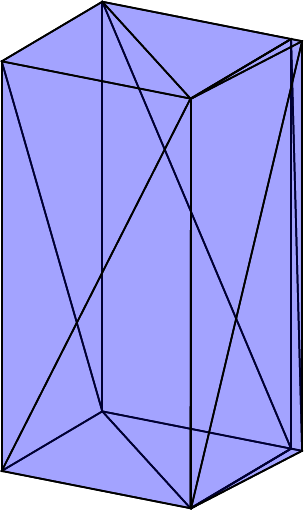}
     \caption{}  \label{fig:tetpoly_mesh_A-b}
   \end{subfigure}
  \caption{Plots of elements from case $A$. (a) Poor-quality tetrahedron and
  (b) Agglomeration of the same tetrahedron with its neighbor to form a polyhedral 
  virtual element.}
           \label{fig:tetpoly_mesh_A}
\end{figure}

Let $\vm{u} \equiv \bigl(u(\vx,t),v(\vx,t),w(\vx,t)\bigr)$ 
be the displacement field and
$\vm{\sigma} \equiv \vm{\sigma}(x,y,z,t)$ be the
Cauchy stress field. For the tapered beam, we solve the following elastodynamic initial/boundary-value problem:
\begin{subequations}
\begin{align}
    \nabla \cdot \vm{\sigma} &
    = \rho \ddot{\vm{u}} 
    \ \ \textrm{in } \Omega \times [0,T_\textrm{max}]\\
    \vm{u}(x,y,z,0) &= 0, \ \ 
    \dot{\vm{u}}(x,y,z,0) = 0, \\
    \vm{u}(0,y,z,t) &= 0,  \ \
    v(4,y,z,t) = 0, \ \ w(4,y,z,t) = 0,  \\
    u(4,y,z,t) &= \begin{cases}
    \left( \dfrac{t}{\tau} \right)^4
    -2  \left( \dfrac{t}{\tau} \right)^3
    +  \left( \dfrac{t}{\tau} \right)^2 & \textrm{for }
    t < \tau \\
    0 & \textrm{otherwise}
    \end{cases},
\end{align}
\end{subequations}
where we choose $\tau = 100 \Delta t_{\textrm{crit}}^{\textrm{VEM}}$ in the numerical simulations.  In addition, we set Poisson's ratio $\nu = 0$ to produce
a two-dimensional solution for the wave propagation problem.

In the numerical computations, an explicit central-difference scheme is used 
for the time integration.\cite{Hughes:2000:FEM} 
We normalize
the displacement in the $x$-direction by $u_0 = 1/16$ and the time by
$T = 4/c_L = 7.7 \times 10^{-4}$, where $c_L = 5188.75$ m/s is the longitudinal wave speed. 
In~\fref{fig:ux_vs_time},
the time history of the normalized displacement in the $x$-direction
at the nodal location $(2,1/2,0)$ is presented. We observe that the virtual
element and finite element solutions are in good agreement. 
In~\fref{fig:ux_vs_time-a}
(case $A$), the critical time step of the VEM is 275 times that
of the critical time step in the FEM.  In~\fref{fig:ux_vs_time-b} (case $B$), 
the speed-up of the virtual
element simulation (453 time steps) over finite elements (20.9 million time steps) is
of ${\cal O}(10^4)$. For this problem, we point out that it is impractical to use the
element-eigenvalue inequality to set the time step, which further reinforces the advantages of VEM.
We emphasize that in both the FEM and VEM meshes there are nearly co-located nodes connected by an edge with \proofcor{an} extremely small length compared to the average element size.  It is known that the stable time step in the FEM is highly dependent on this edge length; however, it is remarkable that the stable time step in the VEM is essentially independent of it.
\begin{figure}
  \centering
  \begin{subfigure}{0.46\textwidth}
  \includegraphics[width=\textwidth]{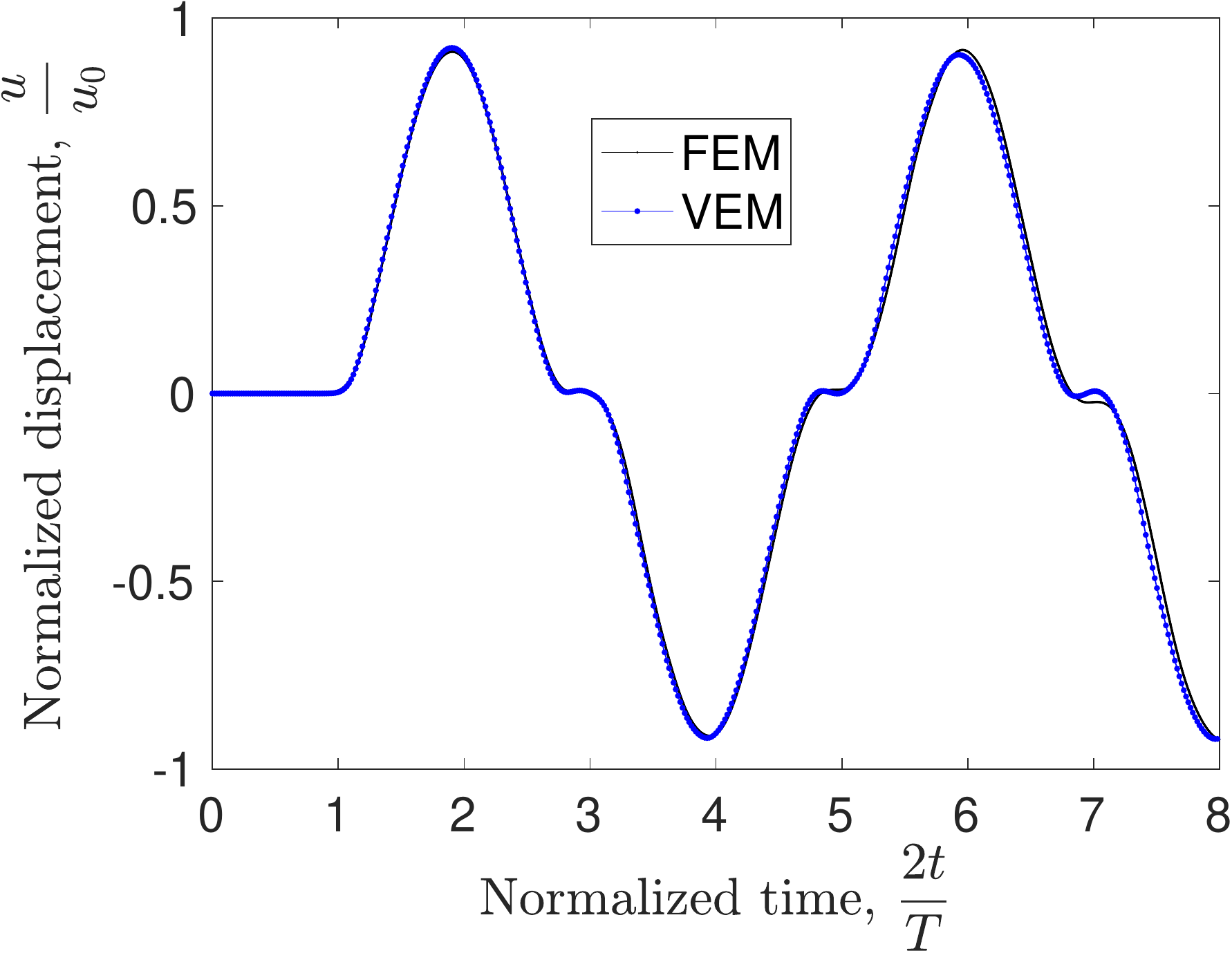}
  \caption{} \label{fig:ux_vs_time-a}
  \end{subfigure} \hspace*{0.1in}
   \begin{subfigure}{0.46\textwidth}
  \includegraphics[width=\textwidth]{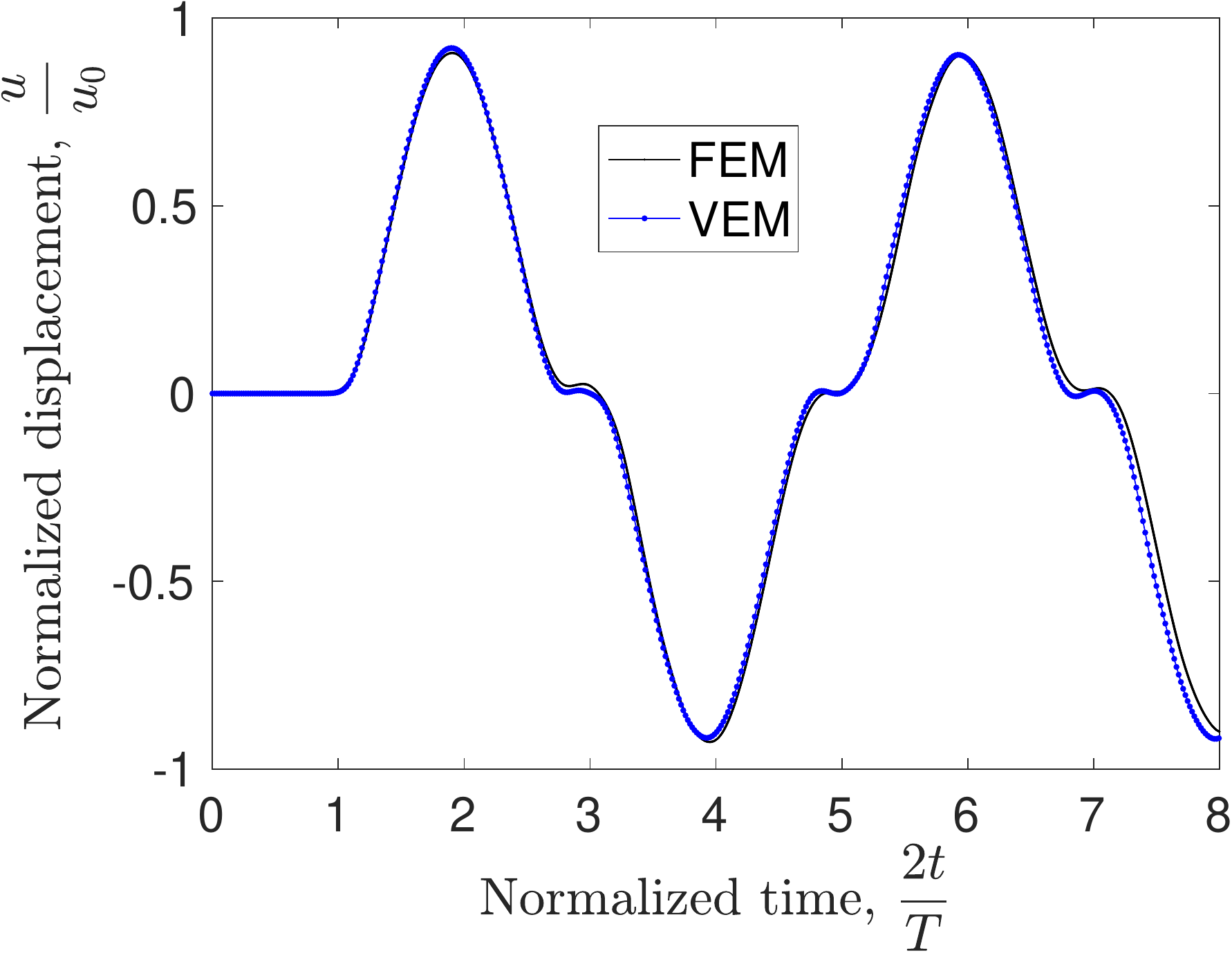}
  \caption{} \label{fig:ux_vs_time-b}
  \end{subfigure}
  \caption{Time history of the normalized displacement in
           the $x$-direction. Comparisons of FE and virtual element solutions for 
           the meshes in case $A$ (Figures~\ref{fig:mesh_AB-a} 
           and~\ref{fig:mesh_AB-b}) are shown in (a).  Solutions for the
           meshes in case $B$ (Figures~\ref{fig:mesh_AB-c} 
           and~\ref{fig:mesh_AB-d}) are shown in (b). In (b), the speed-up of virtual element simulation versus FEM is of ${\cal O}(10^4)$.
           }
           \label{fig:ux_vs_time}
\end{figure}
Finally, we point out that consistent with expectations, instability in the
time history of the 
FE solution is observed when
$\Delta t_{\textrm{FEM}} = 1.000001 \Delta t_{\textrm{crit}}^{\textrm{FEM-g}}$.
}

\section{Conclusions}\label{sec:conclusions}
In this paper, we adopted the first-order
virtual element method
to increase the critical time step for linear 
elastodynamic finite element simulations on tetrahedral meshes that contain poor-quality elements (small faces and dihedral angles close to $0^\circ$ and/or $180^\circ$). We considered tetrahedra with wedge, sliver and spire elements.\cite{Cheng:2000:SE} To this end, we
agglomerated poorly-shaped prismatic and tetrahedral 
elements with their neighbors into polyhedral virtual elements to estimate the critical time step. Lumped mass matrices using row-sum (diagonal scaling) technique for convex (nonconvex) polyhedral elements were utilized.\cite{Park:2019:ONM} The element-eigenvalue inequality~\cite{Fried:1972:BEE,Lin:EET:1989,Hughes:2000:FEM} was used to provide an estimate of the critical time step. For prismatic, wedge and sliver elements, we found that agglomeration with
neighboring tetrahedra into polyhedral elements (with finite volume) produced
critical time steps that were insensitive as the mesh 
parameter $\epsilon \rightarrow 0$. For spire tetrahedral elements, agglomeration produced mixed results:
poor performance resulted when the volume of the polyhedral virtual element approached zero
as $\epsilon \rightarrow 0$, but for an agglomerated polyhedral virtual element with finite volume,
the critical time step was found to be insensitive as 
$\epsilon \rightarrow 0$. Notably, in the favorable instances, the ratio of the critical time step of the polyhedral virtual element to that of the badly-shaped finite element was ${\cal O}(\epsilon^{-1})$, which increases as $\epsilon \rightarrow 0$. This observation suggests that
agglomeration should be done with neighboring elements so that the polyhedral volume does not vanish as $\epsilon \rightarrow 0$. 
\suku{Lastly,
we presented explicit dynamics simulations on a three-dimensional tapered beam
to demonstrate that the increase in critical time step using polyhedral virtual elements leads to significantly faster computer
simulations when compared to tetrahedral finite elements.}

This study has shown the promise of virtual element
technology to enable accurate and reliable linear
elastodynamic finite element simulations on low-quality
tetrahedral meshes. 
Future work will focus on studying the influence of
polyhedral shape quality measures on the eigenspectrum and exploring
robust approaches to perform agglomeration of poor-quality tetrahedra into
polyhedral virtual elements for elastodynamic simulations.

\section*{Acknowledgements}
\ack{
NS acknowledges the research support of Sandia National Laboratories 
\suku{to the University of California at Davis}.
The authors thank Eric Chin for 
providing the code to numerically integrate monomials over polytopes using the homogeneous numerical integration 
method. Helpful discussions with Franco Dassi are also gratefully acknowledged.
Sandia National Laboratories is a multimission laboratory managed and operated by National Technology
and Engineering Solutions of Sandia, LLC, a wholly owned subsidiary of Honeywell International, Inc., for
the U.S. Department of Energy's National Nuclear Security Administration under contract DE-NA-0003525.
This paper describes objective technical results and analysis. Any subjective views or opinions that might
be expressed in the paper do not necessarily represent the views of the U.S. Department of Energy or the
United States Government.
}


\begin{thebibliography}{10}
\providecommand \doibase [0]{http://dx.doi.org/}%

\bibitem{Hormann:2017:GBC}
Hormann K, Sukumar N. \kern-2pt, eds.{\it Generalized Barycentric Coordinates
  in Computer Graphics and Computational Mechanics}.
\newblock Boca Raton: Taylor \& Francis, CRC Press .
\newblock 2017.

\bibitem{Beirao:2013:BPV}
{Beir\~ao da Veiga} L, Brezzi F, Cangiani A, Manzini G, Marini LD, Russo A.
  Basic principles of virtual element methods. {\it Math Models Methods Appl
  Sci} 2013\string; 23\string: 119--214.

\bibitem{Beirao:2014:MFD}
{Beir{\~a}o da Veiga} L, Lipnikov K, Manzini G. {\it The Mimetic Finite
  Difference Method for Elliptic Problems}. 11 of {\it MS\&A -- Modeling,
  Simulation and Applications}.
\newblock Cham: Springer .
\newblock 2014.

\bibitem{Flanagan:1981:AUS}
Flanagan DP, Belytschko T. A uniform strain hexahedron and quadrilateral with
  orthogonal hourglass control. {\it Int J Numer Methods Eng} 1981\string;
  17(5)\string: 679--706.

\bibitem{Cangiani:2015:HSV}
Cangiani A, Manzini G, Russo A, Sukumar N. Hourglass stabilization and the
  virtual element method. {\it Int J Numer Methods Eng} 2015\string;
  102(3--4)\string: 404--436.

\bibitem{Ahmad:2013:EPV}
Ahmad B, Alsaedi A, Brezzi F, Marini LD, Russo A. Equivalent projectors for
  virtual element methods. {\it Comput Math Applications} 2013\string;
  66\string: 376--391.

\bibitem{Beirao:2014:HGV}
{Beir\~{a}o da Veiga} L, Brezzi F, Marini LD, Russo A. The hitchhiker's guide
  to the virtual element method. {\it Math Models Methods Appl Sci}
  2014\string; 24(8)\string: 1541--1573.

\bibitem{Beirao:2017:HOV}
{Beir\~{a}o da Veiga} L, Dassi F, Russo A. High-order {Virtual Element Method}
  on polyhedral meshes. {\it Comput Math Applications} 2017\string; 74\string:
  1110--1122.

\bibitem{Dassi:2018:EHO}
Dassi F, Mascotto L. Exploring high-order three dimensional virtual elements:
  {B}ases and stabilizations. {\it Comput Math Applications} 2018\string;
  75(9)\string: 3379--3401.

\bibitem{Beirao:2013:VEL}
{Beir\~{a}o da Veiga} L, Brezzi F, Marini D. Virtual elements for linear
  elasticity problems. {\it SIAM J Numer Anal} 2013\string; 51(2)\string:
  794-812.

\bibitem{Gain:2014:VEM}
Gain AL, Talischi C, Paulino GH. On the {V}irtual {E}lement {M}ethod for
  three-dimensional linear elasticity problems on arbitrary polyhedral meshes.
  {\it Comput Methods Appl Mech Eng} 2014\string; 282\string: 132--160.

\bibitem{Artioli:2017a:AO2}
Artioli E, {Beir{\~a}o da Veiga} L, Lovadina C, Sacco E. Arbitrary order 2D
  virtual elements for polygonal meshes: part {I}, elastic problem. {\it Comput
  Mech} 2017\string; 60(3)\string: 355--377.

\bibitem{Mengolini:2019:EPV}
Mengolini M, Benedetto MF, Arag{\'o}n AM. An engineering perspective to the
  virtual element method and its interplay with the standard finite element
  method. {\it Comput Methods Appl Mech Eng} 2019\string; 350\string:
  995--1023.

\bibitem{Park:2019:ONM}
Park K, Chi H, Paulino GH. On nonconvex meshes for elastodynamics using virtual
  element methods with explicit time integration. {\it Comput Methods Appl Mech
  Eng} 2019\string; 356\string: 669--684.

\bibitem{Park:2020:NRE}
Park K, Chi H, Paulino GH. Numerical recipes for elastodynamic virtual element
  methods with explicit time integration. {\it Int J Numer Methods Eng}
  2020\string; 121(1)\string: 1--31.

\bibitem{Antonietti:2021:AVE}
Antonietti PF, Manzini G, Mazzieri I, Mourad HM, Verani M. The arbitrary-order
  virtual element method for linear elastodynamics models: convergence,
  stability and dispersion-dissipation analysis. {\it Int J Numer Methods Eng}
  2021\string; 122(4)\string: 934--971.

\bibitem{Cihan:2020:VEF}
Cihan M, Aldakheel F, Hudobivnik B, Wriggers P. Virtual element formulation for
  finite strain elastodynamics. arXiv preprint: 2002.02680;  2020.

\bibitem{Koester:2019:ADS}
Koester JJ, Tupek MR, Mitchell SA. An agile design-to-simulation workflow using
  a new conforming moving least squares method. Tech. Rep. SAND2019-11851,
  Sandia National Laboratories; Albuquerque, NM 87185, USA:   2019.

\bibitem{Bassi:2012:OTF}
Bassi F, Botti L, Colombo A, Petro DAD, Tesini P. On the flexibility of
  agglomeration based physical space discontinuous {Galerkin} discretizations.
  {\it J Comput Phys} 2012\string; 231(1)\string: 45--65.

\bibitem{Cangiani:2014:HPD}
Cangiani A, Georgoulis EH, Houston P. $hp$-version discontinuous {Galerkin}
  methods on polygonal and polyhedral meshes. {\it Math Models Methods Appl
  Sci} 2014\string; 24(10)\string: 2009--2041.

\bibitem{Bishop:2020:PFE}
Bishop JE, Sukumar N. Polyhedral finite elements for nonlinear solid mechanics
  using tetrahedral subdivisons and dual-cell aggregation. {\it Comput Aided
  Geom Des} 2020\string; 77\string: 101812.

\bibitem{Shewchuk:WGL:2002}
Shewchuk JR. What is a good linear finite element? {Interpolation},
  conditioning, anisotropy, and quality measures (preprint). Department of
  Computer Science, University of California, Berkeley, CA 94720, USA;  2002.

\bibitem{Klingner:2008:ITM}
Klingner BM. {\it Improving Tetrahedral Meshes}. PhD thesis. Department of
  Electrical Engineering and Computer Sciences, University of California,
  Berkeley, CA 94720, USA;  2008.

\bibitem{Gillette:2017:SQG}
Gillette A, Rand A. Shape quality for generalized barycentric interpolation. in
  Hormann and Sukumar \cite{Hormann:2017:GBC}ch.~2\string: 23--42.

\bibitem{Attene:2021:BPQ}
Attene M, Biasotti S, Bertoluzza S, et al. Benchmarking the geometrical
  robustness of a {Virtual Element Poisson} solver. {\it Math Comput
  Simulation} 2021\string; 190\string: 1392--1414.

\bibitem{Fried:1972:BEE}
Fried I. Bounds on the extremal eigenvalues of the finite element stiffness and
  mass matrices and their spectral condition numbers. {\it Journal of Sound and
  Vibration} 1972\string; 22(4)\string: 407--418.

\bibitem{Lin:EET:1989}
Lin JJ. An element eigenvalue theorem and its application for stable time
  steps. {\it Comput Methods Appl Mech Eng} 1989\string; 73\string: 283--294.

\bibitem{Lin:BEF:1991}
Lin JJ. Bounds on eigenvalues of finite element systems. {\it Int J Numer
  Methods Eng} 1991\string; 32\string: 957--967.

\bibitem{Benvenuti:2019:EVE}
Benvenuti E, Chiozzi A, Manzini G, Sukumar N. Extended virtual element method
  for the {Laplace} problem with singularities and discontinuities. {\it Comput
  Methods Appl Mech Eng} 2019\string; 356\string: 571--597.

\bibitem{Vacca:2015:VEM}
Vacca G, {Beir{\~a}o da Veiga} L. Virtual element methods for parabolic
  problems on polygonal meshes. {\it Numer Meth Part D E} 2015\string;
  31(6)\string: 2110--2134.

\bibitem{Hughes:2000:FEM}
Hughes TJR. {\it The Finite Element Method: Linear Static and Dynamic Finite
  Element Analysis}.
\newblock Mineola, NY: Dover Publications, Inc. .
\newblock 2000.

\bibitem{Hinton:1976:NML}
Hinton E, Rock T, Zienkiewicz OC. A note on mass lumping and related processes
  in the finite element method. {\it Earthquake Engineering \& Structural
  Dynamics} 1976\string; 4(3)\string: 245--249.

\bibitem{Chin:2015:NIH}
Chin EB, Lasserre JB, Sukumar N. Numerical integration of homogeneous functions
  on convex and nonconvex polygons and polyhedra. {\it Comput Mech}
  2015\string; 56(6)\string: 967--981.

\bibitem{Chin:2020:AEM}
Chin EB, Sukumar N. An efficient method to integrate polynomials over polytopes
  and curved solids. {\it Comput Aided Geom Des} 2020\string; 82\string:
  101914.

\bibitem{Cheng:2000:SE}
Cheng SW, Dey TK, Edelsbrunner H, Facello MA, Teng SH. Sliver exudation. {\it
  Journal of the ACM} 2000\string; 47(5)\string: 883--904.

\end{thebibliography}
\end{document}